\documentclass[12pt]{amsart}
\usepackage{amsmath}
\numberwithin{equation}{section}
\usepackage{amsthm}
\usepackage{amssymb}
\usepackage{amsfonts}
\usepackage{pinlabel}
\usepackage{graphicx}
\usepackage{enumerate}
\usepackage{caption}
\usepackage{subcaption}
\usepackage{xcolor}
\usepackage{color}
\newcommand{\diam}{\mathrm{diam}}
\newcommand{\ra}{{\rightarrow}}
\newcommand{\eproof}{\hfill\rule{2.2mm}{3.0mm}}

\newcommand{\Area}{\mathrm{Area}}

\newcommand{\Vol}{\mathrm{Vol}}
\newcommand{\R}{{\mathbf R}}

\newcommand{\Z}{{\mathbf Z}}

\newcommand{\Q}{{\mathbf Q}}
\newcommand{\N}{{\mathbf N}}
\renewcommand{\H}{{\mathbf H}}

\newcommand{\PSL}{\mathrm{PSL}}

\newcommand{\D}{\noindent}

\newcommand{\grad}{\mathrm{grad}}
\newcommand{\id}{\mathrm{id}}

\renewcommand{\eqref}[1]{(\ref{#1})}
\newcommand{\eat}[1]{}

\newcommand{\bsmall}{\begin{array}[c]{c}}
\newcommand{\esmall}{\end{array}}

\newcommand{\vol}{\mathrm{Vol}}
\newcommand{\inter}{\mathrm{int}}

\newcommand{\inj}{\mathrm{inj}}

\renewcommand{\div}{\mathrm{div}}
\newcommand{\tphi}{\tilde{\phi}}
\newcommand{\tpsi}{\tilde{\psi}}
\newcommand{\tD}{\tilde{D}}
\newcommand{\tS}{\tilde{S}}
\newcommand{\tf}{\tilde{f}}
\newcommand{\tg}{\tilde{g}}

\theoremstyle{plain}
\newtheorem{theo}{Theorem}[section]
\newtheorem{lem}[theo]{Lemma}

\newtheorem{coro}[theo]{Corollary}

\theoremstyle{definition}

\newtheorem{conj}[theo]{Conjecture}
\newtheorem{exam}[theo]{Example}
\newtheorem{Claim}[theo]{Claim}

\newtheorem*{rmk}{Remark}

\newtheorem{fact}[theo]{Fact}

\newcommand{\F}{{\mathcal F}}

\def\>{>_{\sigma}}

\title{The Cheeger Constant, Isoperimetric Problems, and Hyperbolic Surfaces}
\author{Brian Benson}
\address{Department of Mathematics\\
Kansas State University\\
Manhattan, KS 66506}
\email{babenson@ksu.edu}

\begin{document}
\baselineskip 18pt

\begin{abstract}
We give a brief literature review of the isoperimetric problem and discuss its relationship with the Cheeger constant of Riemannian $n$-manifolds. For some non-compact, finite area 2-manifolds, we prove the existence and regularity of subsets whose isoperimetric ratio is equal to the Cheeger constant. To do this, we use results of Hass-Morgan for the isoperimetric problem of these manifolds. We also give an example of a finite area 2-manifold where no such subset exists. Using work of Adams-Morgan, we classify all such subsets of geometrically finite, finite area hyperbolic surfaces where such subsets always exist. From this, we provide an algorithm for finding these sets given information about the topology, length spectrum, and distances between the simple closed geodesics of the surface. Finding such a subset allows one to directly compute the Cheeger constant of the surface. As an application of this work suggested by Agol, we give a test for Selberg's eigenvalue conjecture. We do this by comparing a quantitative improvement of Buser's inequality resulting from works of both Agol and the author to an upper bound on the Cheeger constant of these surfaces, the latter given by Brooks-Zuk. As expected, our test does not contradict Selberg's conjecture.
\end{abstract}
\maketitle

\section{Introduction}
\label{sec:intro}
\setcounter{equation}{0}

The {\bf Cheeger constant} of a finite volume Riemannian $n$-manifold $M$ is given by 
	\begin{equation*}
		h(M):=\inf_D \frac{\Vol_{n-1}(\partial D)}{\Vol_n(D)}
	\end{equation*} 
where $D \subset M$ is a smooth $n$-submanifold with boundary and $0< \Vol_n(D) \leq \Vol_n(M)/2$.
We describe the procedure of using results about the isoperimetric problem on $M$ in order to prove the existence and regularity of subsets $A \subseteq M$ with $0<\Vol_n(A) \leq \Vol_n(M)/2$ so that $$h(M)=\frac{\Vol_{n-1}(\partial A)}{\Vol_n(A)}.$$ We refer to $\partial A$ and $A$ as {\bf $(n-1)$- and $n$-dimensional Cheeger minimizers} respectively. Buser introduced this procedure when he proved that Cheeger minimizers exist for all compact Riemannian manifolds \cite[Remark 3.3, Lemma 3.4]{B82}. Since Buser's explanation of this idea was very brief, we elaborate on it in the setting of compact Riemannian manifolds in Section \ref{sec:isop}. We also discuss criteria for when these results  extend to non-compact, finite volume surfaces.

In section \ref{sec:noncom}, we focus on Riemannian 2-manifolds which we call surfaces. In Example \ref{ex:nonconst2}, we exhibit a non-compact, finite area surface which does not have a Cheeger minimizer. We define a Cheeger sequence $D_k$ to be a sequence of smooth subsurfaces with boundary of $S$ such that $\Area(D_k) \leq \Area(S)/2$ and $\lim_{k\to \infty} h^{\ast}(D_k)=h(S)$ where $h^{\ast}(D_k)$ denotes the isoperimetric ratio of $D_k$. Using work of Hass and Morgan \cite{HM}, we then prove the existence of Cheeger minimizers for finite area, non-compact surfaces with a Cheeger sequence whose boundaries are uniformly bounded in $S$:

\noindent {\bf Theorem \ref{theo:Cheeger}.} {\em Let $S$ be a finite area surface. If there exists a Cheeger sequence $D_k$ such that the boundaries $\partial D_k$ are uniformly bounded in $S$, then $S$ has a 2-dimensional Cheeger minimizer $D$ so that $\partial D$ is an embedded multi-curve of constant curvature.}
	
Adams and Morgan give a complete classification of solutions to the isoperimetric problem in geometrically finite hyperbolic surfaces \cite{AM}. Using their classification and the relationship between the isoperimetric problem and the Cheeger constant, we give a classification of Cheeger minimizers of geometrically finite hyperbolic surfaces in Section \ref{sec:class}. We call sets which have the form of a Cheeger minimizer from our classification {\bf good Cheeger candidates}. We then use Theorem \ref{theo:Cheeger}, to prove the following:

\noindent {\bf Theorem \ref{theo:curvh}.}
   {\em Let $S$ be a geometrically finite hyperbolic surface. Then there exists a 1-dimensional Cheeger minimizer in $S$ with positive length which is the boundary of a good Cheeger candidate of $S$.}

While Cheeger's inequality holds for both compact and non-compact manifolds, Buser gives a different version of his inequality for the case of a non-compact manifold without boundary \cite{B82}. In general, the version of Buser's inequality for closed manifolds does not hold when the manifold is not compact. For closed manifolds, quantitative improvements to Buser's inequality were first considered by Agol \cite{IA}. The author reformulated these improvements in terms of an eigenvalue of a Sturm-Liouville problem. This Sturm-Liouville eigenvalue depends on the identical geometric invariants of $M$ used by Buser's inequality to produce an upper bound on $\lambda_1(M)$ \cite[Theorem 1.4]{BB}.  These results do not depend directly on the compactness of the manifold, instead, they depend on the existance of an $(n-1)$-dimensional Cheeger minimizer. Therefore, Theorem \ref{theo:curvh} allows the application of the quantitative improvements of Buser's inequality to geometrically finite hyperbolic surfaces.

In Section \ref{sec:algorithm}, we use the classification of Cheeger minimizers from Section \ref{sec:class} to give an algorithm for directly computing the Cheeger constant of a geometrically finite hyperbolic surface $S$. To implement this algorithm, we assume that we know the lengths of and distances between all geodesics of $S$ which satisfy a length bound depending on the area of $S$. This is done by finding or approximating a Cheeger minimizer.

Selberg's eigenvalue conjecture suggests that $\lambda_1(S_k) \geq \frac{1}{4}$ for the hyperbolic surfaces $S_k$ which are specific quotients of congruence subgroups. In Section \ref{sec:Selberg}, we give an application of Theorem \ref{theo:Cheeger} and Corollary \ref{coro:Surf}, along with a reformulation of work of Agol \cite{IA} appearing in Benson \cite{BB}, to produce a lower bounds on $h(S_k)$ depending on a lower bound on $\lambda_1(S_k)$. When we assume that $\lambda_1(S_k) \geq \frac{1}{4}$, our lower bound on $h(S_k)$ is consistent with an upper bound on $h(S_k)$ given by Brooks and Zuk \cite{BZ02}, as suggested by Selberg's conjecture.

\section{Background and Motivation}
\label{sec:back}
\setcounter{equation}{0}

Recall that the {\bf Cheeger constant} of a finite volume Riemannian $n$-manifold $M$ is defined to be  
	\begin{equation*}
		h(M):=\inf_D \frac{\Vol_{n-1}(\partial D)}{\Vol_n(D)}
	\end{equation*} 
with $D \subset M$ a smooth $n$-submanifold with boundary and $0< \Vol_n(D) \leq \Vol_n(M)/2$.
We will assume that all manifolds and surfaces are connected, unless otherwise stated. For $f \in C^2(M)$, we denote the Laplacian of $f$ by $\Delta f=-\div \big (\grad (f)\big )$. When $M$ is compact, we denote by $\lambda_1(M)$ the smallest positive value $\lambda \in \R$ so that $\Delta f = \lambda f$ for a function $f \in C^2(M)$. The eigenvalue $\lambda_1(M)$ is called the smallest non-zero eigenvalue of $M$. When $M$ is non-compact, one defines
	$$\lambda (M)=\inf_f \frac{\int_M \|\grad (f) \|^2 \, d\Vol_n}{\int_M f^2 \, d\Vol_n}$$
where $f$ runs over all $C^1(M)$.\footnote{It is more natural to require that $f$ is in the Sobolev space $H^1(M)=W^{1,2}(M)$, as distinct eigenfunctions are linearly independent and form a basis for this Hilbert space.} Since $M$ has finite volume, we require that $\int_M f \, d\Vol_n =0$.

Cheeger gave the initial motivation for the Cheeger constant by showing that 
$$\lambda_1(M)\geq \frac{h(M)^2}{4}$$ \cite{C69}.
 Buser then showed that for $M$ a closed Riemannian $n$-manifold with Ricci curvature bounded below by $ -\delta^2(n-1)$ where $\delta \geq 0$, then $$\lambda_1 (M) \leq 2\delta (n-1) h(M)+10h^2(M)$$ \cite[Theorem 1.2]{B82}. Although we will discuss the Cheeger constant of such manifolds herein, Buser gave an example showing that this inequality does not hold when the manifold has non-empty boundary \cite{B82}. Agol gave a quantitative improvement to Buser's inequality \cite{IA}; see Benson for details \cite[Theorem 1.2]{BB}. Agol's result was motivated by hyperbolic 3-manifolds where his method produced a threefold improvement. In previous work of the author, we showed that Agol's result is equivalent to a statement that $\lambda_1(M)$ is bounded above by the first eigenvalue of a Sturm-Liouville eigenvalue problem, where this eigenvalue problem depends on $h(M)$ as a parameter \cite[Theorem 1.4]{BB}. Specifically, the following eigenvalue comparison holds: 
	\begin{theo} \label{theo:AB} {\bf (Agol \cite{IA}, Benson \cite{BB})}
		Given the assumptions on $M$ in Buser's inequality above, there is an explicit Sturm-Liouville problem depending on $h=h(M)$, $n$, and $\delta$, call it $\omega (h)$ with first eigenvalue denoted by $\lambda_1\big (\omega (h) \big )$. Then
		\begin{equation}\label{eq:AB}
		\lambda_1(M) \leq \lambda_1\big (\omega (h)\big ).
		\end{equation}
	\end{theo}
We describe the precise formulation of $\omega (h)$ for finite area hyperbolic surfaces in Section \ref{sec:Selberg}. While $\lambda_1 \big (\omega (h) \big )$ is not written as an explicit function, one can show that it is a real-valued, continuous, differentiable a.e. function of $h$ for fixed $n$ and $\delta$. This follows from combining the work of Atkinson \cite{FA64}, Everitt, Kwong, and Zettl \cite{EKZ}, Kong, Wu, and Zettl \cite{KWZ}, and M\"oller and Zettl \cite{MZ96}. See Benson for more details \cite[Section 4.4]{BB}.

A key component of the proof of Theorem \ref{theo:AB} is the existence of an $(n-1)$-dimensional Hausdorff measurable set $\Sigma$ dividing $M$ into two components with the following properties:
\begin{enumerate}
	\item The set $\Sigma= \partial A$ for some $n$-dimensional Hausdorff measurable set $A \subset M$ with $\Vol_n(A) \leq \Vol_n(M)/2$.
	\item The sets $\Sigma$ and $A$ achieve the Cheeger constant; that is
		$$h(M)= \frac{\Vol_{n-1} (\Sigma)}{\Vol_n(A)}.$$
\end{enumerate}
When $M$ is closed, such a $\Sigma$ is guaranteed to exist by the combined work of many authors, see Theorem \ref{theo:GMT1} and Lemma \ref{lem:Buser}. The relationship between the isoperimetric problem and the Cheeger constant established by Buser \cite{B82} plays a crucial role in proving this existence. When $M$ is non-compact, such a $\Sigma$ need not exist. Indeed, in Example \ref{ex:nonconst2}, we give an example of a non-compact, finite-volume $2$-manifold where no such $\Sigma$ exists.

\section{The Isoperimetric Problem and the Cheeger Constant}
\label{sec:isop}
\setcounter{equation}{0}

Let $M$ be an $n$-dimensional, Riemannian manifold with $\Vol_n (M) < \infty$.  For a positive $n$-Hausdorff measure subset $A \subseteq M$ having $(n-1)$-Hausdorff measurable boundary $\partial A$, the {\bf isoperimetric ratio} of $A$ is given by
	\begin{equation}\label{eq:hast}
		h^{\ast} (A) = \frac{\Vol_{n-1} (\partial A)}{\Vol_n (A)}.
	\end{equation}
The isoperimetric ratio is well-defined when $A$ is a smooth submanifold with boundary of $M$ or, more generally, when $A$ is an integral current.
We also define the {\bf isoperimetric constant} of a subset $A \subseteq M$ to be
	\begin{equation} 
		\bar{h}(A)=\inf_{D \subseteq A} \frac{\Vol_{n-1} (\partial D)}{\Vol_n (D)},
	\end{equation}
where $D$ ranges over all subsets of $A$ which have smooth boundary.

For $t \in \big ( 0, \Vol_n (M) \big )$, the isoperimetric problem of volume $t$ studies the problem of finding the subset $A$ of $M$ with $\Vol_n (A)=t$ and 
	$$\Vol_{n-1}(\partial A) = \inf \big \{\Vol_{n-1} (\partial B):\Vol_n(B)=t \big \}.$$
We call such an $A$ an {\bf isoperimetric minimizer}. For an in-depth exposition on the isoperimetric problem for Riemannian manifolds, see Ros \cite{Ros}.

The Cheeger constant is related to the isoperimetric problem in that we wish to minimize the isoperimetric ratio among all isoperimetric problems corresponding to volumes $t \in \big (0, \Vol_n(M)/2 \big ]$. We will call a set $D$ in $M$ with $h^{\ast}(D)$ well-defined and $\Vol_n(D) \leq \Vol_n(M)/2$ a {\bf Cheeger candidate}.

As described by Ros \cite{Ros}, when combining the results of Almgren \cite{FA76}, Gr\"uter \cite{G87}, and Gonzalez, Massari, and Tamanini \cite{GMT}, one obtains the {\bf fundamental existence and regularity theorem} for isoperimetric minimizers of Riemannian manifolds.
\begin{theo}\label{theo:GMT1} {\bf (Fundamental Existence and Regularity)}
	When $M$ is smooth, compact, possibly with boundary, and $t\in \big (0, \Vol_n(M) \big )$, there exists a compact region $A \subset M$ such that $\partial A$ minimizes the $(n-1)$-volume among all regions of $n$-volume $t$. Further, for any such minimizer, $\partial A$ is a smooth embedded hypersurface of constant mean curvature, except for a singular set of Hausdorff dimension at most $n-8$. Finally, if $\partial M \cap \partial A \neq \emptyset$, then $\partial M$ and $\partial A$ meet orthogonally.
\end{theo}
At least two additional contributions to Theorem \ref{theo:GMT1} should be mentioned. First, Bombieri simplified portions of Almgren's work relevant to the isoperimetric problem \cite{EB82}. Second, Morgan clarified how the results of the other authors, formulated for subsets of Euclidean space, can be adapted to Riemannian manifolds \cite{FM03}.

As suggested by Theorem \ref{theo:GMT1} and smooth approximation techniques, it is equivalent to minimize the isoperimetric ratio over smooth submanifolds with boundary or the more general collection of rectifiable currents to find the Cheeger constant.
In the same way that we have defined isoperimetric minimizers, we can define minimizers for the Cheeger constant, which we refer to as Cheeger minimizers for brevity. Specifically, a {\bf Cheeger minimizer} is a subset $A$ of $M$ such that $\Vol_n(A)\leq \Vol_n(M)/2$ and $h^{\ast}(A) = h(M)$. When $M$ is compact, Buser showed that there exist volume(s) $t \in \big (0,\Vol_n(M)/2 \big ]$ so that Cheeger minimizers coincide with solutions to the isoperimetric problem of $n$-volume $t$ in $M$ \cite[Remark 3.3, Lemma 3.4]{B82}. Therefore, the existence and regularity of a Cheeger minimizer follows from Theorem \ref{theo:GMT1}.
Specifically, Buser proved:
	\begin{lem}\label{lem:Buser} {\bf (Buser \cite[Lemma 3.4]{B82})}
		Let $M$ be an n-dimensional, compact Riemannian manifold, possibly with boundary. Then there exists $t^{\ast} \in \big (0, \Vol_n (M)/2 \big ]$ and a sequence of submanifolds $A_k \subset M$ with smooth boundary such that $\Vol_n (A_k) =t^{\ast}$ for all $k$ and
			$$h(M)=\lim_{k \ra \infty} \frac{\Vol_{n-1} (\partial A_k)}{\Vol_n (A_k)}.$$
	\end{lem}
	
Combining Theorem \ref{theo:GMT1} and Lemma \ref{lem:Buser}, we can immediately conclude that Cheeger minimizers exists for all compact manifolds. The following example shows that such a result is not always possible for non-compact, finite volume manifolds:

\begin{exam}\label{ex:nonconst2} We will construct a surface $S_0$ such that $h(S_0)=0$. Denote the disk of Euclidean radius $r$ centered at the origin of $\R^2$ by $D_r$. Now consider a region $R$ given by $\R^2-D_1^{\circ}$, where $D_1^{\circ}$ is the interior of $D_1$. The metric on $R$ will be of the form $ds^2=dr^2+\frac{1}{r^4} d\theta^2$ in the Euclidean polar coordinates $(r, \theta)$. This metric is complete since every Cauchy sequence is contained in a coordinate annulus of the form $\big \{(r, \theta):0<r_1<r<r_2<\infty\big \}$. For $r>1$, it follows that $h^{\ast} \left (D_r^{\complement}\right ) = \frac{1}{r} \to 0$ as $r \to \infty$. 

We can then extend the metric on $R$ to the complement of $R$, which we denote $R^{\complement}$, by letting the metric on $R^{\complement}$ be the inversion of the metric $ds^2$ on $R$ with respect to the variable $r$. We then smooth the metric in a small tubular neighborhood of the circle $r=1$ in polar coordinates. This gives a metric on all of $\R^2-\{(0,0)\}$ and we denote the resulting surface  $S_0$. By our computations above, we can conclude that $h(S_0)=0$ since each $D_r^{\complement}$ is a Cheeger candidate for $r >1$. This follows from the construction of our metric, since the area on the interior of the circle of Euclidean radius 1 is equal to the area on the exterior of the circle of Euclidean radius 1. Finally, the following computation shows that $S_0$ has finite area:
	$$\Area(S_0)=2\int_0^{2\pi} \int_1^{\infty} \frac{1}{r^2} \, dr \, d\theta= 4\pi.$$
\end{exam}

\begin{rmk} This is also an example of a finite area surface with a Cheeger constant of 0. Buser proved that this is not possible for compact manifolds \cite{B82}. Note that it is not necessary for the areas of a sequence of sets with isoperimetric ratios converging to the Cheeger constant to go to infinity for the Cheeger constant to be zero.
\end{rmk}

In response to Example \ref{ex:nonconst2}, we consider the properties or criteria one can specify on a non-compact finite volume manifold to guarantee the existence of a Cheeger minimizer. As a first step towards answering this question, we will give two related criteria, one is a special case of the other (see Lemma \ref{lem:BoundCom}), which guarantee the existence of a Cheeger minimizer in dimension 2; see Section \ref{sec:noncom}. We say that a sequence of subsets $D_k$ of $M$ is {\bf uniformly bounded} if there exists a compact subset $K$ of $M$ such that $D_k \subseteq K$ for all $k \in \N$. Further, we will say that a sequence of Cheeger candidates $A_k \subset M$ is a {\bf Cheeger sequence} in $M$ if $\Vol_n(A_k) \leq \Vol_n(M)/2$ and $\lim_{k \ra \infty} h^{\ast}(A_k)=h(M)$. We will show that a version Lemma \ref{lem:Buser} holds for non-compact finite volume $n$-manifolds under the additional assumption that the boundaries of the Cheeger sequence are uniformly bounded:

\begin{lem}{\bf (Equal Volume)}\label{lem:EqVol}
	Let $M$ be a non-compact, finite volume Riemannian $n$-manifold. Assume there exists a Cheeger sequence $A_k$ such that the $\partial A_k$ are uniformly bounded. Then there exists a Cheeger sequence of $n$-submanifolds with boundary $D_k \subset M$ such that $\Vol_n(D_k)=\liminf_{i \to \infty} \Vol_n(A_i)=:t$ for all $k \in \N$. Consequently, a Cheeger minimizer is given by a solution to the isoperimetric problem of volume $t$.
\end{lem}

In order to prove Lemma \ref{lem:EqVol}, we show that the volumes of a Cheeger sequence cannot vanish whenever the boundaries of the sequence are uniformly bounded. A version of the following result is given by Buser in his proof of Lemma \ref{lem:Buser} for compact manifolds: 

\pagebreak
\begin{lem}\label{lem:Bnonvanish}{\bf (Non-Vanishing Volume)}
	Let $M$ be a non-compact, finite-volume Riemannian $n$-manifold. If $A_k$ is a Cheeger sequence in $M$ and the boundaries $\partial A_k$ of the sequence are uniformly bounded in $M$, then $$\liminf_{k} \Vol_n(A_k)>0.$$
\end{lem}

Both Lemma \ref{lem:EqVol} and Lemma \ref{lem:Bnonvanish} will be proved in Section \ref{sec:Proofs}.

\subsection{The Cheeger Constant of Flat Tori and the Klein Bottle}
In this section, we will give a short example to illustrate the method of applying Buser's Lemma \ref{lem:Buser} to examples of surfaces for which the isoperimetric problem has been solved for all areas. As a result, we can directly compute the Cheeger constant of these surfaces using Howards' solutions to the isoperimetric problem for flat tori and the Klein bottle \cite{H92}; see also Hutchings, Howards, and Morgan \cite[Section 7]{HHM}.

\begin{exam}\label{ex:Tori}
Consider a flat torus and a flat Klein bottle whose metrics are inherited by a rectangle of parallel side lengths $a$ and $b$, with $a \leq b$, in the Euclidean plane with opposite sides identified to one another in the usual ways. The solution of the isoperimetric problem for all areas on flat tori and the Klein bottle are due to Howards \cite{H92}; see also Howards, Hutchings, and Morgan \cite[Section 7]{HHM}. Specifically, for areas $t \in (0,a^2/\pi)$, solutions are metric disks and, for areas $t \in \big (a^2/\pi, \Area (S)-a^2/\pi \big )$, solutions are a band (possibly M\"obius in the case of the Klein bottle) with geodesic boundary components of length $a$. Since each of these surfaces is compact, Buser's Lemma \ref{lem:Buser} implies that a Cheeger minimizer for each is given by a particular solution to the isoperimetric problem. Since $\Area(S) \geq a^2$, the isoperimetric ratio as a function of $t$ is monotone decreasing for $t\in \big (0, \Area(S)/2\big ]$. Therefore, we may conclude that a band bounded by two closed geodesics of length $a$ enclosing an area equal to $\Area(S)/2$ is a Cheeger minimizer. As a result, the Cheeger constant of these surfaces is equal to $\frac{2a}{ab/2}=4/b$.
\end{exam}

\subsection{Preliminary Results for Isoperimetric Ratios}\label{sec:Prem}
\begin{figure}[ht!]
\labellist
\small \hair 2pt
\pinlabel $B'$ at 10 475
\pinlabel $A$ at 175 375
\pinlabel $B$ at 200 70
\pinlabel $A'$ at 720 85
\pinlabel $A$ at 660 330
\pinlabel (a) at 200 0
\pinlabel (b) at 770 0
\endlabellist
\centering
\includegraphics[width=\textwidth]{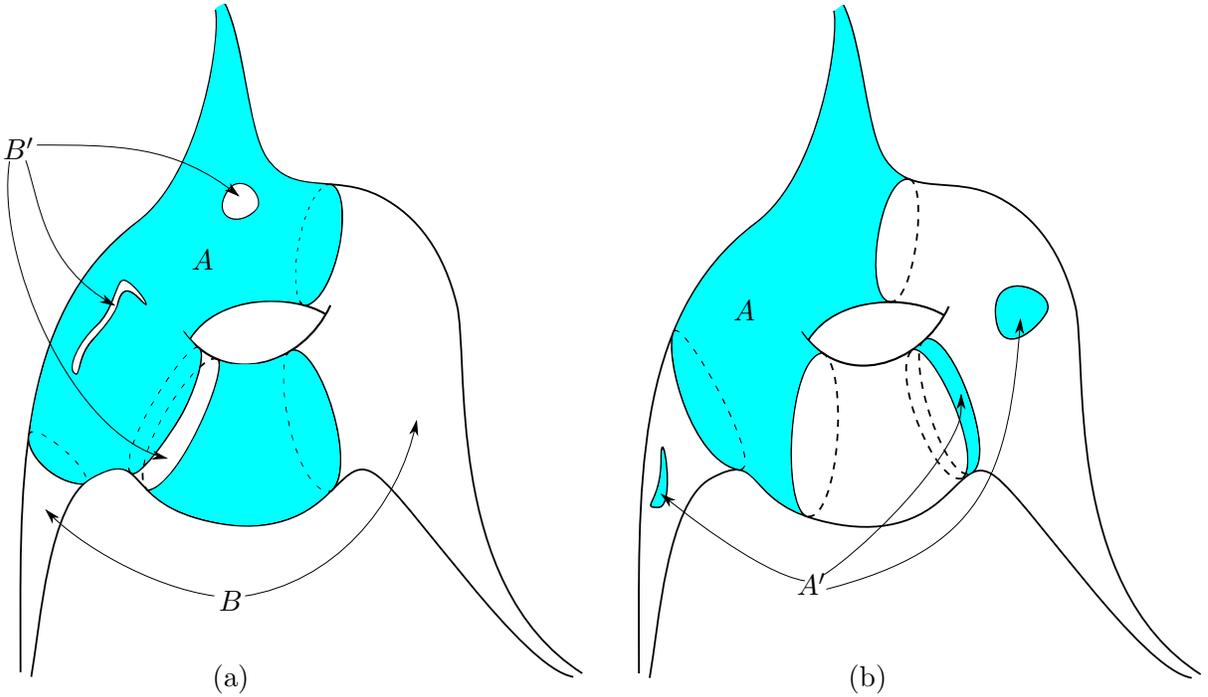}
\caption{(a) An example of $A$, $B$, and $B'$ from Lemma \ref{lem:CIL} where $A$ is shaded. (b) An example of $A$ and $A'$ from Lemma \ref{lem:CEL} where $A$ is shaded.}
\label{fig:CILCEL}
\end{figure}

We will prove two lemmas which allow us to simplify our assumptions about Cheeger sequences. Examples of the sets described in these lemmas are illustrated in Figure \ref{fig:CILCEL}. Suppose that $\Sigma$ is an $(n-1)$-dimensional Cheeger candidate and the two components of $M-\Sigma$ are $A$ and $B$ with $\Vol_n(A) \leq \Vol_n(B)$. In this situation, the Component Inclusion Lemma \ref{lem:CIL} gives an isoperimetric estimate where we want to include an $n$-dimensional component to $A$ from $B$ to decrease the $(n-1)$-volume of $\Sigma$:
\begin{lem}\label{lem:CIL} {\bf (Component Inclusion Lemma)} Let $A$ be a Cheeger candidate for $M$.  Let $B'$ be a union of connected components of $M - A$ and $B$ be all the other connected components of $M-A$. If $h^{\ast} (B') > h^{\ast} (A)$, then at least one of the following is true:
	\begin{enumerate}
		\item $h^{\ast}(A \cup B') < h^{\ast} (A)$ and $\Vol_n(A \cup B') \leq \frac{1}{2} \Vol_n (M)$,
		\item $h^{\ast} (B) < h^{\ast} (A)$ and $\Vol_n (B) \leq \frac{1}{2} \Vol_n (M)$.
	\end{enumerate}
\end{lem}
 The Component Exclusion Lemma \ref{lem:CEL} is more straight-forward to prove and handles the estimate when we wish to remove an $n$-dimensional component from $A$ and include it in $B$ to decrease the $(n-1)$-volume of $\Sigma$:
 \begin{lem}\label{lem:CEL} {\bf (Component Exclusion Lemma)} Let $A$ be a Cheeger candidate of $M$. Let $A'$ be a union of connected components of $A$. If $h^{\ast} (A')>h^{\ast} (A)$, then $h^{\ast}(A - A') < h^{\ast}(A)$.
\end{lem}
 
In proving these lemmas, we will use two straightforward facts.
	\begin{fact}\label{fact1}
		If $B$ and $B'$ are disjoint, then 
			$$h^{\ast} (B \cup B') \geq \min \big \{ h^{\ast}(B), h^{\ast} (B')\big \}$$
		with equality exactly when $h^{\ast}(B)=h^{\ast}(B')$.
	\end{fact}
	\begin{fact}\label{fact2}
		If $A$ is a Cheeger candidate, then $h^{\ast} (A) \geq h^{\ast}(M-A)$.
	\end{fact}

\noindent {\bf Proof of Lemma \ref{lem:CIL}.}  Since $M=(A \cup B') \cup B$, and $A,B',B$ are all pairwise disjoint, at least one of the following cases is true: either $\Vol_n (A \cup B') \leq \frac{1}{2} \Vol_n(M)$ or $\Vol_n (B) \leq \frac{1}{2}\Vol_n(M)$. First, we have the following comparison of isoperimetric ratios:
	\begin{equation}\label{eq:AcupB'}
		h^{\ast}(A \cup B') = \frac{\Vol_{n-1}(\partial( A \cup B'))}{\Vol_n (A \cup B')} = \frac{\Vol_{n-1} (\partial A) - \Vol_{n-1}( \partial B')}{\Vol_n (A \cup B')} < \frac{\Vol_{n-1} (\partial A)}{\Vol_n (A)}
		= h^{\ast} (A).
	\end{equation}
In addition, by Facts \ref{fact1} and \ref{fact2}, 
	$$h^{\ast}(A) \geq h^{\ast}(M-A)= h^{\ast}(B \sqcup B') \geq \min \big \{ h^{\ast} (B), h^{\ast} (B') \big \}.$$ 
Since $h^{\ast} (B') > h^{\ast} (A)$, we must have $h^{\ast} (A)  > h^{\ast} (B)$ as a consequence of Fact \ref{fact1} when $h^{\ast}(B) \neq h^{\ast}(B')$. Specifically, the strict inequality follows by Fact \ref{fact1}, because $h^{\ast}(B') > h^{\ast} (B)$, so $h^{\ast} (B \sqcup B') > h^{\ast} (B)$.
\eproof

\noindent {\bf Proof of Lemma \ref{lem:CEL}.} We have $$h^{\ast}(A)=h^{\ast} \big ( A' \sqcup (A-A') \big ) \geq \min \big \{ h^{\ast}(A'), h^{\ast}(A-A')\big \}$$ by Fact \ref{fact1}. Further, $h^{\ast}(A') > h^{\ast}(A)$ by hypothesis, so $h^{\ast}(A-A') < h^{\ast}(A')$ and we conclude that $h^{\ast}(A)>h^{\ast}(A-A')$ by the equality condition of Fact \ref{fact1}.
\eproof

\subsection{Proofs of Lemmas \ref{lem:EqVol} and \ref{lem:Bnonvanish}} \label{sec:Proofs} Using the results from Section \ref{sec:Prem}, we will now prove Lemma \ref{lem:EqVol} and Lemma \ref{lem:Bnonvanish}. Since the proof of Lemma \ref{lem:EqVol} depends on Lemma \ref{lem:Bnonvanish}, we will first prove Lemma \ref{lem:Bnonvanish}. The portion of the proof for the case where $A_k$ is contained in a compact subset of $M$ is due to Buser \cite{B82}.

\noindent {\bf Proof of Lemma \ref{lem:Bnonvanish}.} Let $K$ be a compact submanifold with boundary of $M$ so that $\partial A_k$ is properly contained in $K$. By replacing $K$ with a superset of itself if necessary, we may assume that each connected component of $M-K$ is non-compact. For each $k \in \N$ where $A_k$ is not completely contained in $K$, we know that $A_k$ must contain at least one connected component of $K^{\complement}$. Since each component of $K^{\complement}$ has positive volume, we can conclude that $\Vol_n(A_k)$ is bounded below by the least volume connected component of $K^{\complement}$.

To complete the proof, we consider the case where $A_k$ is properly contained in $K$ and follow arguments given by Buser for closed manifolds \cite[Proof of Lemma 3.4]{B82}. Using a triangulation of $K$, we can write $K=T_1 \cup \cdots \cup T_l$ where each $T_i$ is mapped homeomorphically onto the Euclidean ball of radius 1, which we denote $B_1^n$, by a fixed quasi-isometry $\Phi_i$ and where $\inter(T_i) \cap \inter(T_j) = \emptyset$ whenever $i\neq j$. If $A_k$ has sufficiently small $n$-volume, then $\Vol_n\big (\Phi_i(A_k \cap T_i) \big ) \leq \Vol_n(B_1^n)$ for $i = 1,\ldots , l$. By the classical isoperimetric inequality for $B_1^n \subset \R^n$, for each $i$ where $\Vol_n(A_k \cap T_i) >0$:
	$$\frac{\Vol_{n-1} \Big (\Phi_i \big ( \inter(T_i)\cap \partial A_k \big )  \Big )}{\Vol_n\big ( \Phi_i(T_i \cap A_k \big )^{1+\frac{1}{n}}}\geq C_1^i$$
for a dimensional constant $C_1^i$. Note that, for each $k$, $\Vol_n(A_k \cap T_i) >0$ must occur for at least one $i$.

It follows that for each $i=1,\ldots, l$ with $\Vol_n(T_i \cap A_k)>0$:
	\begin{align*}
	\frac{\Vol_{n-1}\big (\inter (T_i) \cap \partial A_k \big )}{\Vol_n(T_i \cap A_k)}&\geq C^i_2 \Vol_n(A_k\cap T_i)^{-\frac{1}{n}}\\ 
	&\geq C^i_2 \big (\Vol_n(A_k) \big )^{-\frac{1}{n}}
	\end{align*}
where $C_2^i$ depends on the length distortion of $\Phi_i$, but is independent of $A_k$. Then we have that 
	\begin{align*}
	h^{\ast}(A_k)&=\frac{\sum_i \Vol_{n-1}\big (\inter (T_i)\cap \partial A_k \big )}{\sum_i \Vol_n(T_i \cap A_k)}\\
	&\geq \min_i \frac{\Vol_{n-1}\big (\inter(T_i) \cap \partial A_k \big )}{\Vol_n(T_i \cap A_k)}\\ & =C\big (\Vol_n(A_k) \big )^{-\frac{1}{n}}
	\end{align*}
where $C=\min_i C^2_i$ and the minima are taken over each $i=1, \ldots, l$ so that $\Vol_n(T_i \cap A_k)>0$. It follows that 
	\begin{equation*}
		\infty > h(M)=\liminf_{k} h^{\ast}(A_k)\geq C \liminf_{k} \big ( \Vol_n(A_k)\big )^{-\frac{1}{n}},
	\end{equation*}
and so $\liminf_k \Vol_n(A_k)>0$.
\eproof

We again use arguments provided by Buser for closed manifolds, this time to prove Lemma \ref{lem:EqVol}. Specifically, our proof of Lemma \ref{lem:EqVol} follows Buser's idea of attaching or removing small metric balls to the sets in the Cheeger sequence and show this gives a Cheeger sequence where all the sets have the same volume:

\noindent {\bf Proof of Lemma \ref{lem:EqVol}.} Let $A_k$ be a Cheeger minimizing sequence. By the Non-Vanishing Lemma \ref{lem:Bnonvanish}, define $t:=\liminf_k h^{\ast}(A_k)>0$. For any small $\epsilon>0$, we can divide the interval $[t, \Vol_n(M)/2]$ into closed subintervals of length smaller than $\epsilon$ (which will have pairwise intersection at the endpoints). At least one of these subintervals must contains the quantities $\Vol_n (A_k)$ for infinitely many $k \in \N$. Further reducing the $A_k$ to such a subsequence, we may assume that $\Vol_n (A_k)$ differ by less than $\epsilon$ from one another for all $k\in \N$. To avoid double counting volume, we will assume that the radius of a ball of $n$-volume $\epsilon$ is smaller than the injectivity radius of each point in the set $\bigcup_{p \in \partial A_k} B_{\epsilon}(p)$ for all $k \in \N$. Such an $\epsilon$ exists since the $\partial A_k$ are uniformly bounded and so it is possible to find a compact submanifold with boundary of $M$ containing $\bigcup_{p \in \partial A_k}B_{\epsilon} (p)$. By further restricting to a subsequence if necessary, we may assume that the quantities $\Vol_n(A_k)$ are monotone.

Take a small $r_0>0$, $p_k \in A_k$, and $q_k \in A_k^{\complement}$ so that $B_{r_0}(q_k) \subset A_k^{\complement}$, $B_{r_0}(p_k) \subset A_k$, and there exists a unique distance minimizing geodesic $\gamma_k:[0,1] \to S$ with $\gamma_k (0)=p_k$ and $\gamma_k(1)=q_k$. Since $\partial A_k$ is contained in a compact subset of $S$ and $\Vol_{n-1}(\partial A_k)$ is bounded, we can find an $r_0$ where these properties hold for all $k \in \N$. Further, $\inf_{s,k} \Vol_n \big (B_{r_0}(\gamma_k(s))\big )>0$ since $\bigcup_{k \in \N} \bigcup_{s \in [0,1]} B_{r_0} \big ( \gamma_k(s)\big )$ are contained in a compact subset of $S$, so take $\epsilon$ such that $0<\epsilon<\inf_{s,k} \Vol_n \big (B_{r_0}(\gamma_k(s))\big )$. Then $\Vol_n \big ( B_{r_0}(\gamma_k(s))\big )>\epsilon>\Vol_n(A_1)-\Vol_n(A_k):=d_k$ for all $k \in \N$ and all $s \in [0,1]$.

Define $f(s):=\Vol_n \big (A_k \cup B_{r_0}(\gamma (s)) \big )$ which is a continuous function in $s$ with $f(0)=0 <t+d_k$ and $f(1)>t+d_k$. By the intermediate value theorem, there exists an $s_k \in (0,1)$ such that $f(s_k)=t+d_k$. By taking $\epsilon$ small enough compared to $\inf_{s,k} \vol_n \big ( B_{r_0}(\gamma_k (s))\big )$ and $\gamma_k$ so that $\gamma_k(s)$ intersects $\partial A_k$ for exactly one $s \in [0,1]$, we can take $p_k :=\gamma_k(s_k) \in A_k$.

Now $g(r):=\Vol_n\big (A_k \cup B_r(p_k) \big )$ is a continuous function in $r$ with $g(0)=\Area (A_k) <t$ and $g(r_0)=t+\epsilon>t$. Again by the intermediate value theorem, there exists $r_k \in (0, r_0)$ such that $g(r_k)=t$. Define $D_k$ to be a smooth area preserving perturbation of $A_k \cup B_{r_k}(p_k)$ so that $D_k=A_k \cup B_{r_k}(p_k)$ off of a tiny neighborhood of $\partial A_k \cap \partial B_{r_k}(p_k)$.

Since $\Vol_n (A_k) \to t$ as $k \to \infty$, we have that $r_k \to \infty$. Therefore, the metric balls $B_{r_k}(p_k)$ are close to Euclidean balls for all large $k \in \N$. So then there exists a constant $\delta$ so that 
	\begin{equation}\label{eq:Bbound}
	\Vol_{n-1}\big ( \partial B_{r_k}(p_k)\big ) \leq \delta \alpha r^{n-1}_k
	\end{equation}
for all large $k$, where $\alpha$ is the $(n-1)$-volume of the boundary of and $\omega$ is the $n$-volume of the unit Euclidean $n$-ball. As a result, we have that 
\begin{equation*}
\lim_{k \to \infty} h^{\ast}(D_k) \leq \lim_{k \to \infty} \frac{\Vol_{n-1}(\partial A_k) + \Vol_{n-1} (\partial B_{r_k}(p_k))}{\Vol_n(A_k)}.
\end{equation*}
Applying Equation \ref{eq:Bbound} to the previous equation, we have 
\begin{equation*}
 \lim_{k \to \infty} h^{\ast}(D_k) \leq \lim_{k \to \infty}  \frac{\Vol_{n-1}(\partial A_k) + \delta\alpha r_k^{n-1}}{\Vol_n(A_k)}.
 \end{equation*}
 Since $\liminf_{k \to \infty} \Vol_n(A_k)>0$, it follows that 
\begin{equation*}
\lim_{k \to \infty} h^{\ast}(D_k)= \lim_{k \to \infty} \frac{\Vol_{n-1}(\partial A_k)}{\Vol_n (A_k)}=h(M)
\end{equation*}
and $h^{\ast}(D_k)\geq h(M)$ for all $k \in \N$ by the definition of the Cheeger constant since $\Vol_n(D_k)=t\leq \Vol_n(M)/2$. Therefore, we conclude that $\lim_{k \to \infty}h^{\ast}(D_k) = h(M)$. A similar argument can be used for the case when $\Vol_n(A_k)>t$.

\begin{figure}[ht!]
\labellist
\small \hair 2pt
\pinlabel $A_k$ at 40 120
\pinlabel $A_k^{\complement}$ at 50 80
\pinlabel $p_k$ at 130 215
\pinlabel $q_k$ at 196 52
\pinlabel $\gamma_k$ at 177 85
\pinlabel {$\gamma_k (s_k)$} at 163 150
\pinlabel {$B_{r_k}\big ( \gamma_k (s_k)\big )$} at 210 150
\pinlabel {$B_{r_0}(p_k)$} at 200 215
\pinlabel {$B_{r_0}(q_k)$} at 260 80
\endlabellist
\centering
\includegraphics[width=5in]{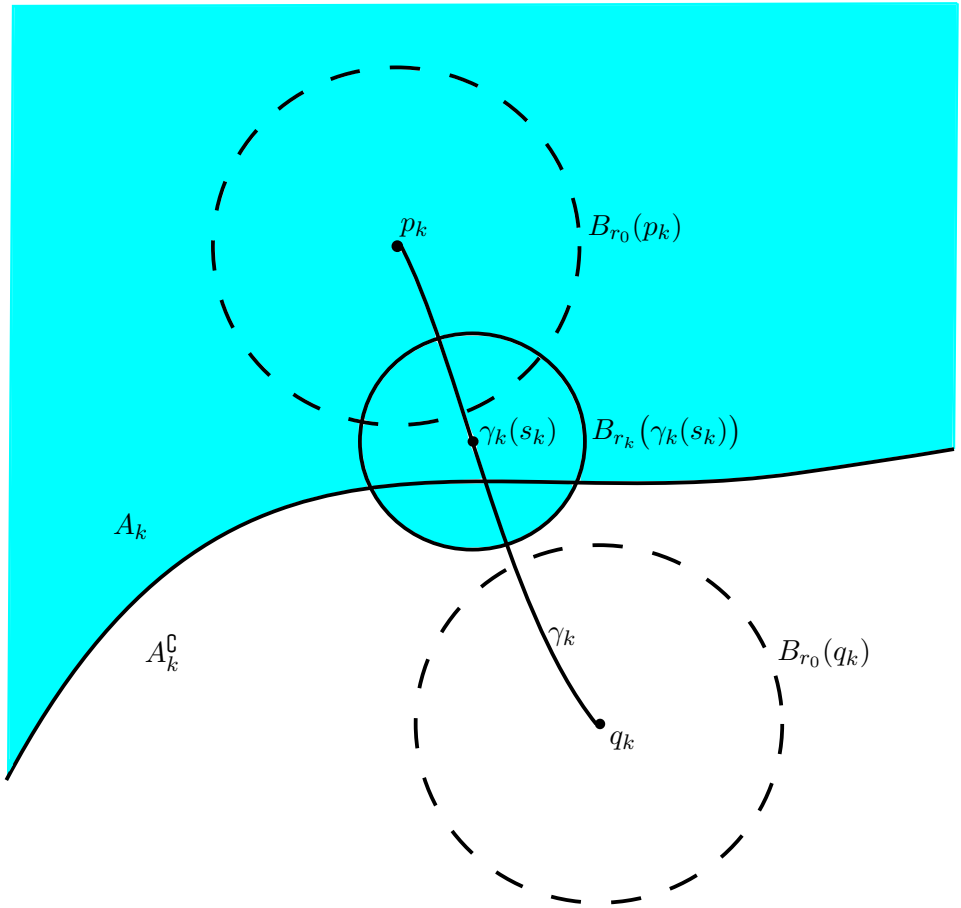}
\caption{A ball of radius $r_k$ is attached to $A_k$ to create the shaded region $D_k=A_k \cup B_{r_k}\big ( \gamma_k(s_k)\big )$. }
\label{fig:Attach}
\end{figure}
\eproof

\section{Non-Compact Finite Area Surfaces}
\label{sec:noncom}
\setcounter{equation}{0}

We will study the Cheeger constant and isoperimetric problem of connected Riemannian 2-manifolds which we call {\bf surfaces}. To simplify notation, we let $\ell := \Vol_1$ and $\Area := \Vol_2$. Using techniques accessible to a wide audience, Hass and Morgan prove the following theorem for surfaces related to Theorem \ref{theo:GMT1}:
	\begin{theo}\label{theo:GMT2} {\bf (Hass and Morgan \cite[Theorem 3.4]{HM})}
		Let $S$ be a smooth closed surface with $t \in \big ( 0, \Area(S) \big )$. Then there is a bounded subsurface with boundary $A$ of $S$ with $\Area (A) =t$ whose boundary $\partial A$ has least length among all subsurfaces with area equal to $t$. Further, its boundary $\partial A$ consists of embedded curves all having the same constant curvature.
	\end{theo}

We consider extending Theorem \ref{theo:GMT2} to some non-compact, finite area surfaces. For such a surface $S$, we call a sequence $D_i$ of submanifolds of $S$ with boundary an {\bf isoperimetric sequence of area} $t$ if $\Area(D_i)=t$ and $\ell (\partial D_i)$ converges to the minimum boundary length for all subsurfaces with boundary of $S$ having area $t$. Using theorems of Hass and Morgan, we prove the following result about solutions to the isoperimetric problem which allows us to prove the existence and regularity of some Cheeger minimizers in Section \ref{sec:ExReg}:
\begin{theo}\label{theo:HM}
	Let $S$ be a finite area, non-compact surface. Suppose that for some $t \in \big ( 0, \Area (S)\big )$, there exists an isoperimetric sequence $D_i$ of submanifolds of $S$ with boundary having area $t$ such that $\partial D_i$ is uniformly bounded in $S$. Then there exists a subset $D$ of $S$ with $\Area (D)=t$ such that $D$ is an isoperimetric minimizer of area $t$ in $S$ and $\partial D$ is an embedded multi-curve of constant curvature in $S$.
\end{theo}

The uniformly bounded assumption of this theorem can be removed by appealing directly to the methods and techniques of Hass and Morgan. See Section \ref{sec:HM} for a description. As we have seen in Example \ref{ex:nonconst2}, the existence of a Cheeger minimizer cannot be guaranteed for all non-compact, finite area surfaces. However, Theorem \ref{theo:HM} and Lemma \ref{lem:EqVol} tells us that a Cheeger minimizer must exist whenever we can find a Cheeger minimizing sequence with uniformly bounded boundary curves:

\pagebreak
\begin{theo}\label{theo:Cheeger}
		Let $S$ be a finite area surface. If there exists a Cheeger minimizing sequence $D_k$ such that the boundaries $\partial D_k$ are uniformly bounded in $S$, then $S$ has a 2-dimensional Cheeger minimizer $D$ so that $\partial D$ is an embedded multi-curve of constant curvature.
	\end{theo}
	
In Section \ref{sec:IsopCom}, we give one criterion for the isoperimetric constants of non-compact submanifolds of boundary of $S$ which guarantees that $\partial D_k$ are uniformly bounded. This is proved in Lemma \ref{lem:BoundCom} in Section \ref{sec:IsopCom}. We find this criterion convenient for proving the existence and regularity of Cheeger minimizer in finite area hyperbolic surfaces addressed in Section \ref{sec:class}.
	
\subsection{Techniques of Hass and Morgan}\label{sec:HM}

A stronger version of Theorem \ref{theo:HM} is possible without the assumption that the boundaries of a minimizing sequence are uniformly bounded by arguing directly using methods of Hass and Morgan \cite{HM}. While Hass and Morgan do not state such a result in their work, it can be proved directly using techniques provided by them in their proof of Theorem \ref{theo:GMT2}. Specifically, one can prove the following result which we still attribute to Hass and Morgan:
	\begin{theo} \label{theo:HM2} {\bf (Hass and Morgan \cite{HM})}
		Let $S$ be a smooth, finite area, non-compact surface and $t \in \big (0, \Area(S) \big )$. Then there exists a subset $D$ of $S$ with $\Area (D)=$ such that $D$ is an isoperimetric minimizer of area $t$ in $S$ and $\partial D$ is an embedded multi-curve of constant curvature in $S$.
	\end{theo}
In what follows, we will describe how to apply the Hass-Morgan techniques to prove this result.
Note that for any $t \in \big ( 0,\Area(S)\big )$, at least one connected component of the boundary of an isoperimetric sequence must be uniformly bounded. Hass and Morgan apply the Arzela-Ascoli Theorem to prove existence and we may still apply this theorem to the aforementioned uniformly bounded connected components of the boundary of the sequence. Specifically, one first minimizes over sets with fixed area and whose boundary curves share a fixed combination of homotopy types. Then one applies local convexity and first variation arguments to prove regularity, specifically that boundary curves are embedded and have constant curvature except for finitely many segments of multiplicity two. The local convexity argument assures that the limit curves, whose existence are guaranteed by Arzela-Ascoli, do not intersect themselves in such a way that the first variation argument cannot be applied. Using this local convexity argument in both directions along each segment of multiplicity two, one proves that these segments of multiplicity two are geodesics. Now one must minimize over regions bounded by curves of unrestricted homotopy types. A key observation is that one can improve the isoperimetric ratio of a region bounded by curves having segments of multiplicity two by removing such segments. This increases the number of boundary components of the region and their homotopy types. This idea is depicted in Figure \ref{fig:ImagBound}. By the proof that multiplicity one segments of minimizers of fixed homotopy types are embedded and have constant curvature, it suffices to prove that the number of boundary curves must be bounded for a region to solve the isoperimetric problem in $S$. In order to bound the number of components of a minimizer for the isoperimetric problem (without homotopy constraints), one can use the classical isoperimetric inequality, an argument that we also use in the proof of Theorem \ref{theo:Cheeger} that we give herein and originally appeared in Hass and Morgan \cite[pages 192-193]{HM}. As observed by Howards, Hutchings, and Morgan, since non-compact subsurfaces with boundary of the surface have finite area, any sequence of boundary curves traveling out to infinity bound an area going to zero towards infinity and, therefore, may be discarded \cite[page 434]{HHM}. It follows that one may ignore sequences of connected components of the boundary which are not uniformly bounded.

\subsection{Existence and Regularity of Isoperimetric Minimizers}\label{sec:ExReg}
Recall that, by Example \ref{ex:nonconst2}, it is possible for there to be finite area surfaces which do not have Cheeger minimizers. To prove Theorem \ref{theo:Cheeger} about the existence and regularity of Cheeger minimizers of finite area surfaces given a boundedness condition on Cheeger sequences, we find that Theorem \ref{theo:HM2} is more general than is necessary. Therefore, in order to keep the literature as complete and self-contained as possible, we  use the specific statements of theorems from the Hass and Morgan paper \cite{HM} about compact surfaces with boundary and show how to extend them to non-compact, finite area surfaces with Cheeger sequences having uniformly bounded boundary curves. We accomplish this by applying Lemma \ref{lem:EqVol} taking a Cheeger sequence which is also an isoperimetric sequence, then refining each set in the sequence by minimizing its isoperimetric ratio over its fixed topological types of its boundary curves in a subset of $S$. By showing that the number of connected components of boundary curves of sets in this refined sequence is bounded, we will show that a minimizer is achieved for some finite index in the refined sequence. 

For a smooth embedding $f:D \to S$, define the isoperimetric ratio of the embedding to be 
	$$h^{\ast}(f):=\frac{\ell \big (f(\partial D) \big )}{\Area \big ( f(D)\big )}.$$
For a smooth curve $\gamma:[0,L] \to S$ and for a point $t_0 \in [0,L]$, denote $\gamma'(t_0):=\gamma_{\ast} \left (\left. \frac{d}{dt}\right |_{t=t_0} \right ) \in T_{\gamma(t_0)} S$ where $\gamma_{\ast}$ is the push-forward of $\gamma$. Recall that the mapping length of $\gamma$ coincides with the Hausdorff measure when $\gamma$ is injective:
	$$\ell(\gamma) =\int_0^L |\gamma'(t)|\, dt.$$
The {\bf uniform norm} of $\gamma$ is given by $\|\gamma\|_{\infty}=\sup\{|\gamma'(t)|:t \in [0,L]\}$. We will now prove Theorem \ref{theo:HM}.

\noindent {\bf Proof of Theorem \ref{theo:HM}.} The existence of an isoperimetric sequence $D_i$ of area $t$ with $\partial D_i \subset K$ for some compact $K \subset S$ is given by assumption. Take $\epsilon>0$ with $\epsilon \ll t$ and consider a compact subsurface with boundary $\tS$ of $S$ which properly contains $K$ as a subset and such that $\Area(S-\tS)<\epsilon$. Denote the isoperimetric ratio of subsets of and functions to $S$ and $\tS$ by $h^{\ast}_S$ and $h^{\ast}_{\tS}$ respectively.

For $\tD_i:=D_i\cap \tS$, the boundaries $\partial \tD_i$ consist of two mutually exclusive types of connected components: components in $\partial D_i $ and components in $\partial \tS$. It is possible that no components of this second type exist. Since $K \subset \tS$, we have $\partial D_i \subseteq \partial \tD_i$.

For connected components $F$ of $\tS^{\complement}$, either $F \subset D_i$ or $F \subset D_i^{\complement}$ for each $i \in \N$. Denote by $I_i$ the set of these components $F$ such that $F \subset D_i$. By passing to a subsequence, we may assume that $I_i=I_j$ for all $i,j \in \N$ and we refer to these sets as $I$ without ambiguity.
	Now we will consider minimizing the isoperimetric ratio over a set $\F_i$ of the closure of smooth embeddings $\tf:\tD_i \to \tilde{S}$ such that 
		\begin{enumerate}
			\item For every connected component $E \subseteq \partial \tS\cap \tD_i$, we have that $\left. \tf \right |_E$ is the identity map. Recall that each connected component of $\partial D_i$ is a subset of $\tS$, so each connected component of $\partial \tS$ is either a subset of $\tD_i$ or a subset of $\tD_i^{\complement}$.
			\item For every connected component $C \subseteq \partial D_i \subset \partial \tD_i$, the image $\tf(C)$ is homotopic to $C$ in $\tS$ when both are parametrized as curves.
			\item $\Area \big ( \tf (\tD_i) \big )=\Area (\tD_i )=\Area (D_i)-\sum_{F \in I} \Area (F) \geq t-\epsilon$.
		\end{enumerate}

	Hass and Morgan prove that there exists a minimizer $\tilde{g}_i \in \F_i$ such that $\tilde{g}_i(\partial \tD_i) \cap \inter (\tilde{S})$ is an embedded multi-curve of constant curvature, except for finitely many geodesic arcs or isolated points where it has multiplicity two \cite[Theorem 3.1, Remark 3.3, Remark 3.5]{HM}. The boundary of every isoperimetric minimizer has constant curvature measured with respect to a consistent choice of an outward pointing normal vector field. Finally, local convexity is in the same direction along connected components of the images of the multi-curve boundary.\footnote{More precisely, the stated result of Hass and Morgan is for the existence and regularity of a map from $\partial \tD_i$ to $\tS$ which is identical to $\left. \tilde{g}_i\right |_{\partial \tD_i}$ for the map we denote $\tilde{g}_i$. We use the statement above to simplify notation and because we do not need to work with the interior of $\tg_i(\tD_i)$ directly.}
	
For each $\phi \in \F_i$ such that $\phi$ is not an embedding, define $h^{\ast}_{\tS}(\phi):=\lim_{k \to \infty} h^{\ast}_{\tS}(\phi_k)$ where $\phi_k$ is a sequence of smooth embeddings in $\F_i$ which converge to $\phi$ in the uniform norm. For any $\tD \subset \tS$, we define $\partial D$ so that $\partial \tS \cap \tD \subseteq \partial \tD$. As a result, we have $\left. \partial D \right |_{\tS} \subseteq \partial \tD$ and so the length of $\partial D$ in $S$ is less than or equal to the length of $\partial \tD$ in $\tS$. It is also straight-forward to see that $\Area (D) \geq \Area (\tD)$. As a result of these inequalities, we write $h^{\ast}_S$ and $h^{\ast}_{\tS}$ to denote the isoperimetric ratios in $S$ and $\tS$ respectively.
	
	Each function $\tf \in \F_i$ can be extended to a map $f:D_i \to S$ by 
		\begin{equation}\label{eq:extend1}
		f(s)=\left \{ \begin{array}{ll}\tf(s),& s \in \tD_i\\
							s,& s \in D_i-\tD_i \end{array}\right.
		\end{equation}
Conditions (1) and (3) imply that $\Area \big (f(D_i) \big )=t$ for every $f \in \F_i$.
Further, for any $\tf \in \F_i$  with extension $f:D_i \to S$, we have
	\begin{equation}\label{eq:extendisop}h_S^{\ast}(f) = \frac{\ell \big ( \tf(\partial \tD_i)\big )
	-\ell \big ( \tf(\partial \tS \cap \tD_i)\big )}
	{\Area\big (\tf(\tD_i) \big )+\Area (D_i-\tD_i)}
	=\frac{\ell \big ( \tf(\partial \tD_i)\big )
	-\ell (\partial \tilde{S} \cap \tD_i)}{t}.\end{equation}
Since the term $\ell (\partial \tS \cap D_i)$ is independent of $\tf$ and $h^{\ast}_{\tS}(\tf)=\ell \big ( \tf(\partial \tD_i)\big )/t$ for every $\phi, \psi \in \F_i$, we can conclude that $h_S^{\ast} (\phi) \leq h_S^{\ast}(\psi)$ if and only if $h^{\ast}_{\tS}(\tphi) \leq h^{\ast}_{\tS}(\tpsi)$. Since we can find a minimizer $\tf_i \in \F_i$ for the sequence $D_i$ with $h^{\ast}_{\tS}(\tf_i) \leq h^{\ast}_{\tS}(\tD_i)$, we have that $h^{\ast}_S(f_i) \leq h^{\ast}_S(D_i)$.

\begin{figure}[ht!]
\labellist
\small \hair 2pt
\pinlabel $\tf_i(\partial \tD_i)$ at 77 125
\endlabellist
\centering
%\vspace{-.5in}
\includegraphics[width=5in]{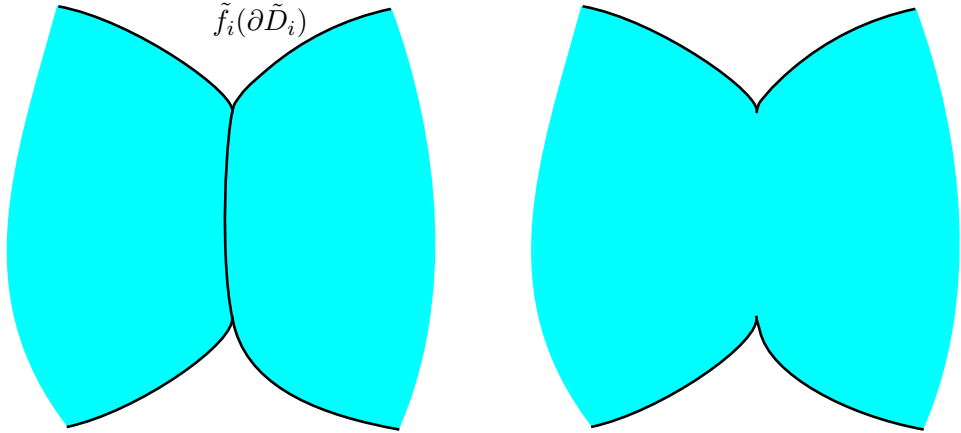}
%\vspace{-1.5in}
\caption{Removing the interior of an arc segment of multiplicity two from the image of $\tf_i(\tD_i)$ reduces length of the boundary while keeping the area of the enclosed region constant.}
\label{fig:ImagBound}
\end{figure}

By construction, the number of connected components of $\partial \tD_i$ is greater than or equal to the number of connected components of $\tf_i(\partial \tD_i)$ and is greater than or equal to the number of connected components of $f_i(\partial D_i)$. Following an argument of Hass and Morgan, we will prove that the number of components of $\tf_i(\partial \tD_i)$, call it $b_i$, is bounded above for all $i \in \N$ \cite[pages 192-193]{HM}. Proceed by contradiction, assuming that the number $b_i$ of connected components of $\tf_i(\partial \tD_i)$ goes to infinity as $i \to \infty$. When $b_i$ is large, then $\ell (D_i) \geq \ell \big (f_i(\partial D_i)\big )$ is close to the minimum, implying that $\ell \big ( \tf_i(\partial \tD_i\big )$ is close to the minimum by Equation \ref{eq:extendisop}. Multiplicity two arc segments of the boundary of $\tf_i(\partial \tD_i)$ can be removed by changing the topology of $D_i$. Specifically, removing the interior of a multiplicity two arc of positive length of the image of $\tf_i(\partial \tD_i)$, as in Figure \ref{fig:ImagBound}, decreases the length of the boundary of the region while leaving the area of the enclosed region unchanged. Therefore, the total length of multiplicity two arc segments of $\tf_i(\partial \tD_i)$ must be small. 

By restricting to a subsequence of $D_i$ if necessary, we may assume that $h^{\ast}(D_i)$ is monotone decreasing and so $\ell \big (\tf_i(\partial \tD_i)\big ) \leq t \cdot h^{\ast}(\tD_1)$ for all $i \in \N$. Since we have assumed that $b_i \to \infty$, a connected component of $\tf_i(\partial \tD_i)$ of multiplicity one must be short and, by the theorem of Hass and Morgan, have the same large constant curvature (in absolute value). So, either $\tf_i(\tD_i)$ or $\tf_i(\tD_i)^{\complement}$ has many connected components which are small disks in $S$. The isoperimetric inequality implies that the isoperimetric ratio of a very small disk increases as the radius of the disk decreases. Since the area enclosed by $\tilde{f}_i(\partial \tilde{D}_i)$ is constant for all $i \in \N$, the lengths $\ell \big (\tf_i(\partial \tD_i) \big )$ must be increasing as $i \to \infty$, a contradiction. So we conclude that $b_i$ has an upper bound $b \in \R$.

Recall that we have $\ell \big (\tf_i(\partial \tD_i)\big ) \leq t \cdot h^{\ast}(\tD_1)$ for all $i \in \N$ by restricting to a subsequence of the $D_i$ if necessary. The genus of each $\tD_i$ is bounded by the genus of $\tS$. Further, $\tS$ is compact and contains finitely many homotopy classes of curves with length less than or equal to $t\cdot h^{\ast} (\tD_1)$ \cite[Remark 1.13 (b), pg. 10]{Grom}. So then there exists only finitely many subsurfaces with boundary $\tD$ of $\tS$ having area $t$ and fewer than $b$ boundary components up to homotopies of connected components of $\partial \tD$. Therefore, for some $j \in \N$, the homotopy classes of the connected components of $\partial \tD_j$ must be in one-to-one correspondence with the homotopy classes of a subsequence $\partial \tD_{i_k}$ (for infinitely many $i_k$). By construction of the $\F_i$, the image of the minimizer $\tf_j \in \F_j$ must have $h^{\ast}_{\tS}(\tf_j)=h^{\ast}_{\tS}(\tf_{i_k})$ for all $i_k$. It follows that $h^{\ast}_S(f_j)=h_S^{\ast}(f_{i_k})\leq h_S^{\ast}(D_{i_k})$ by the extension of these maps on $\tS$ to all of $S$. Since $h_{\tS}^{\ast}(\tD_{i_k})$ and $h^{\ast}_S(D_{i_k})$ will converge to the minimum in $\tS$ and $S$ respectively, we may conclude that the images of $f_j$ and $\tf_j$ are isoperimetric minimizers of $S$ and $\tS$ respectively. Hass and Morgan prove that any isoperimetric minimizer of $\tS$ must have a boundary whose intersection with $\inter(\tS)$ consists of embedded curves all having the same constant curvature \cite[Theorem 3.4, Remark 3.5]{HM}.

To complete the argument, we need to show that $\tf_j(\partial \tD_j)$ consists of embedded curves all having the same constant curvature on all of $\tS$ instead of only on $\inter(\tS)$. To see this, take a smooth submanifold with boundary $\tS^0$ of $S$ with $\tS \subset \inter(\tS^0)$ and define $\tD_j^0:=D_j\cap \tS^0$. Extend $\tf_j:\tD_j \to \tS$ to $\tf^0_j:\tD_j^0 \to \tS^0$ via the relation: 
\begin{equation}\label{eq:extend}
		\tf^0_j(s):=\left \{ \begin{array}{ll}
						\tf_j(s),& s \in \tD_i\\
						s,& s \in \tD_i^0-\tD_i
						\end{array}\right.
	\end{equation}
Define $\F_j^0$ to be the closure (in the uniform norm) of smooth embeddings of the form $\tf^0:\tD^0 \to \tS^0$ in the same way we defined $\F_i$ for smooth embeddings of the form $\tf:\tD_i \to \tS$. Then, because $\tf_j \in \F_j$ extends to an isoperimetric minimizer $f_j:D_j \to S$, then $\tf_j^0= \left. f_j \right |_{\tS^0}$ must minimize the isoperimetric ratio in $\tS^0$ over all functions in the set $\F_j^0$. If not, then there exists $\tg_j^0 \in \F_j^0$ with $h^{\ast}_{\tS^0}(\tg_j^0)<h^{\ast}_{\tS^0} (\tf_j^0)$ which can be extended to a function $g_j:D_j \to S$ via a relation similar to the one given in Equation \ref{eq:extend1}. Further, formulas for $h^{\ast}_S(g_j)$ and $h^{\ast}_S(f_j)$, can be given using the same idea which resulted in Equation \ref{eq:extendisop}. Specifically, we can conclude that the order of $h^{\ast}_S(g_j)$ and $h^{\ast}_S(f_j)$ depend only on the order of $\ell \big (\tg^0_j(\tD_j^0)\big )$ and $\ell \big (\tf_j^0(\tD_j^0)\big )$ respectively. Therefore, we have $h^{\ast}_S(g_j)<h^{\ast}_S(f_j)$, a contradiction. Since $\tf_j^0(\partial \tD_j^0)=f_j(\partial D_j) \subset \inter(\tS_0)$, it follows that $f_j(\partial D_j)$ is an embedded multi-curve of constant curvature by Theorems of Hass and Morgan \cite[Theorem 3.1, Remark 3.3, Theorem 3.4, Remark 3.5]{HM}.
\eproof

\subsection{Existence and Regularity of Cheeger Minimizers} \label{sec:Cheeger}
As is the case with compact manifolds, when they exist, the existence and regularity of Cheeger minimizers of non-compact, finite area surfaces can be obtained using the existence and regularity of isoperimetric minimizers of these surfaces. To see this, we will now complete the proof of Theorem \ref{theo:Cheeger}.

\noindent {\bf Proof of Theorem \ref{theo:Cheeger}.} Lemma \ref{lem:EqVol} tells us that whenever we can find a Cheeger minimizing sequence with uniformly bounded boundary curves, we can find a Cheeger sequence where all of these sets in the sequence have the same area (and the boundary curves remain bounded). Therefore, a Cheeger minimizer must coincide with an isoperimetric minimizer for some $t \in \big (0, \Area(S)/2 \big ]$. The conclusion then follows directly from Theorem \ref{theo:HM}.
\eproof

\subsection{Lower Bounds on the Cheeger Constant}
Upper bounds on the Cheeger constant of $S$ can be obtained by taking the isoperimetric ratio of any Cheeger candidate in $S$. On the other hand, obtaining non-zero lower bounds on the Cheeger constant is less straightforward. We will prove the following lemma which provides such bounds for surfaces and later apply it to the cusps of hyperbolic surfaces addressed in Section \ref{sec:class}.

\begin{lem}\label{lem:Stokes} Let $S$ be a surface. Let $D \subset S$ have finite area and be contained in a neighborhood in $S$ isometric to a disk, possibly with finitely many punctures, with Euclidean polar coordinates $(r, \theta)$ and metric of the form $$ds^2=f_1(r, \theta)^2 dr^2+f_2(r, \theta)^2 d\theta^2.$$ Let $F=F(r,\theta)$ be a function such that 
	\begin{enumerate}
		\item $\frac{\partial F}{\partial r} = f_1(r,\theta)f_2(r,\theta)$,
		\item $|F(r,\theta)|\leq a f_2(r,\theta)$, for some $a \in \R^+$.
	\end{enumerate}
Then $h^{\ast} (D)\geq \frac{1}{a}$.
\end{lem}

\noindent {\bf Proof of Lemma \ref{lem:Stokes}.} First consider the case where $D$ does not have punctures.  Let $\alpha=F(r, \theta)\, d\theta$, so then $d\alpha = \frac{\partial F}{\partial r}\, dr \, d\theta = f_1(r, \theta)f_2(r, \theta) \, dr\, d\theta$, so we have
	\begin{align*}
		\Area(D)&= \int_D d\alpha\\
			&\leq \int_{\partial D} |F(r,\theta)| \, d\theta &&\text{by Stokes' Theorem}\\
			& \leq a \int_{\partial D} f_2(r, \theta) \, d\theta &&\text{by (2)}\\
			& \leq a \int_{\partial D} ds\\
			& = a \cdot \ell (\partial D).
	\end{align*}
Therefore, we can conclude that $h^{\ast} (D) \geq 1/a$ when $D$ does not have punctures.

For the case where $D$ has finitely many punctures, Stokes' Theorem can only be applied to approximations of $D$ arising from taking a compact subset $D'$ of $D$ with boundary. This boundary arises from deleting a disks of arbitrarily small area containing each puncture and we denote the union of these disks by $N$. Letting $D'$ be $D$ minus the interior of $N$. This means that $\partial D'=\partial D \cup \partial N$. So Stokes' Theorem applies to $D'$ and the argument above gives us an estimate on $h^{\ast} (D')$.

Since $D$ has finite area, the areas of $N$ and the lengths of $\partial N$ vanish as one shrinks $N$. Therefore, we take a monotone sequence $D'_k \supsetneq D'_{k+1} \supsetneq \cdots \supsetneq D$ so that
$\Area(D'_k) \nearrow \Area(D)$ and $\ell(\partial D'_k) \searrow \ell(\partial D)$ as $k \ra \infty$. Then we can conclude that $h^{\ast} (D) \geq 1/a$ when $D$ has punctures.
\eproof

\subsection{Cheeger Sequences with Uniformly Bounded Boundary Curves}\label{sec:IsopCom}
Let $S$ be a finite area surface. We consider decompositions of $S$ of the form $S=N \cup S_C$ so that $N$ and $S_C$ are subsurfaces with boundary of $S$ with $N$ non-compact and $S_C$ compact, and $N\cap S_C=\partial N =\partial S_C$ a smooth multicurve. When $\overline{h}(N) > h(S)$ for such a decomposition of $S$, we say that $S$ is {\bf isoperimetrically compact}.

\begin{lem}\label{lem:BoundCom}
	If $S$ be an isoperimetrically compact surface, then there exists a Cheeger sequence $D_k$ such that $\partial D_k$ is uniformly bounded in $S$.
\end{lem}

\noindent {\bf Proof of Lemma \ref{lem:BoundCom}.} Consider a Cheeger minimizing sequence $\tilde{A}_k$ in $S$ where each $\tilde{A}_k$ is connected by Fact \ref{fact1}. We may assume that $h^{\ast}(\tilde{A}_k) < \overline{h}(N)$ where $N$ is the non-compact subsurface with boundary of $S$ given by the isoperimetric compactness of $S$. Therefore, no $\tilde{A}_k$ can be properly contained in $N$.

Let $F_k$ be the union of connected components of $S-\partial \tilde{A}_k$ properly contained in $N$. By our previous observation, note that $F_k \subset A_k^{\complement}$. Further, when $F_k$ is non-empty, we have that $$h^{\ast}(F_k) >h^{\ast} (\tilde{A}_k)>\max \left \{ h^{\ast}(\tilde{A}_k \cup F_k), h^{\ast} \big ( (\tilde{A}_k \cup F_k)^{\complement}\big )\right \}$$ by the Component Inclusion Lemma \ref{lem:CIL}. Let $A_k:=\tilde{A}_k \cup F_k$ when $\Area (\tilde{A}_k \cup F_k) \leq \Area(S)/2$ and $A_k:=(\tilde{A}_k \cup F_k)^{\complement}$ otherwise. It follows that for ever component $C$ of $\partial A_k$, we have $C \cap S_C \neq \emptyset$.

By restricting to a subsequence if necessary, we may assume that $h^{\ast}(A_k)$ is monotone decreasing in $k$. Then for all $k \in \N$, we have that 
	$$\frac{\ell (\partial A_k)}{\Area(S)/2} \leq \frac{\ell (\partial A_k)}{\Area (A_k)}=h^{\ast}(A_k) \leq h^{\ast}(A_1)$$ and so $\ell (\partial A_k) \leq h^{\ast}(A_1)\Area (S)/2$. So defining $L:=h^{\ast}(A_1)\Area (S)/2$, we have that $\ell (\partial A_k)\leq L$ for all $k \in \N$. Let $d:S \to \R$ be the signed distance function from points in $S$ to $\partial S_C$ which is positive for $s\in N$ (i.e. $d^{-1}[0, \infty)=N$). Denote the diameter of $S_C$ by $\diam(S_C)$ and let $K:=d^{-1}[-\diam(S_C),L]$, which is a compact subset of $S$. Then we conclude that $\partial A_k \subset K$ since $\ell (\partial A_k) \leq L$ and $C \cap A_k \neq \emptyset$ for every connected component $C$ of $\partial A_k$.
\eproof

	\begin{coro}\label{coro:Cheeger}
		Let $S$ be an isoperimetrically compact finite area surface. Then $S$ has a 2-dimensional Cheeger miminizer $D$ such that $\partial D$ is an embedded multi-curve of constant curvature.	
	\end{coro}
	
	\noindent {\bf Proof of Corollary \ref{coro:Cheeger}.} This follows directly from Theorem \ref{theo:Cheeger} and Lemma \ref{lem:BoundCom}.\eproof

\begin{rmk} By using a compact surface where a Cheeger minimizer is known, one can sometimes construct examples of non-compact surfaces where Cheeger minimizers exist. Take a compact surface $S_C$ where a Cheeger minimizer $A \subset S_C$ is known. Find a metric disk $D_R(x)$ with $x$ in the interior of $A^{\complement}$ and $R$ small enough so that $D_R(x) \subset A^{\complement}$. Suppose one can replace the metric in $D_R(x)$ with a metric in polar coordinates of the form $ds^2= f_1(r,\theta)^2 dr^2+f_2(r, \theta)^2 d\theta^2$ from Lemma \ref{lem:Stokes} having at least one puncture and so that $a<h^{\ast}(A)$ where $a$ is the quantity related to $f_1,f_2$ in the lemma. Then the resulting surface $S'$ has $h(S')\leq h^{\ast}(A)$, is isoperimetrically compact by Lemma \ref{lem:Stokes}, and has a Cheeger minimizer by Corollary \ref{coro:Cheeger}.
\end{rmk}

\section{Cheeger Minimizers of Hyperbolic Surfaces}
\label{sec:class}
\setcounter{equation}{0}

In this section, we will give several results which provide a framework for the direct computation of the Cheeger constant for geometrically finite, finite area hyperbolic surfaces. By hyperbolic surface, we mean a connected Riemannian 2-manifold with constant section curvature equal to $-1$. If $S$ is a hyperbolic surface, then $S$ is {\bf geometrically finite} if $S=\Gamma \backslash \H^2$ for $\Gamma$ a geometrically finite, torsion-free Fuchsian group. A discrete group $\Gamma$ of isometries of $\H^2$ is {\bf geometrically finite} if it has a fundamental polyhedron $F \subset \H^2$ with finitely many sides, which is convex, and for each side $S$ of $F$, there exists $g \in \Gamma$ such that $S=F \cap gF$. In dimension two, $\Gamma$ geometrically finite is equivalent to $\Gamma$ finitely generated, see Borthwick for more details \cite[Theorem 2.10]{B07}. For more details on Fuchsian groups, see Katok \cite{K92}.

To simplify our presentation, we define a {\bf horocusp neighborhood} in $S$ to be a neighborhood of an end of $S$ which is covered by a horoball in $\H^2$. Adams and Morgan provide a more specific characterization for the isoperimetric problem for hyperbolic surfaces \cite{AM}:
	\begin{theo}\label{theo:GMT3} {\bf (Adams and Morgan \cite[Theorem 2.2]{AM})}
		Let $S$ be a connected, geometrically finite hyperbolic surface. For a given $t \in \big (0, \Area (S) \big )$, a collection of embedded rectifiable curves 				bounding a region $A$ of area $t$ which minimizes $\partial A$ consists of regions of the following four types:
			\begin{enumerate}
				\item a metric ball,
				\item a horocusp neighborhood,
				\item an annulus bounded by two neighboring curves a constant distance from the unique geodesic of the isotopy class containing both curves, 
				\item or regions bounded by one or more geodesics or single neighboring curves, one in each isotopy class of each geodesic, all of which are a single distance $s$ of from their isotopic geodesic, with $s$ consistently oriented into or out of the region.
	\end{enumerate}
		Further, the following inequality holds 
			\begin{equation}\label{eq:length}\ell (\partial A) \leq \sqrt{\Area(A)^2+4\pi \Area (A)}\end{equation}
		with equality in the case of a circle bounding a metric ball. If $S$ has at least one cusp, cases (1) and (3) do not occur and 				$\ell (\partial A) \leq \Area (A)$ with equality from horocyles. Finally, if $\Area (A) <\pi$ and $S$ has cusps, then a minimizer consists of an arbitrary collection of horoball neighborhoods of cusps with boundary having total length of $\Area (A)$.
	\end{theo}

We will call a subset of a hyperbolic surface $S$ which has one of the following forms a {\bf good Cheeger candidate}:
	\begin{enumerate}
		\item metric disks with area equal to $\Area(S)/2$, only in cases where $S$ is compact,
		\item annuli (and their complements) with area equal to $\Area (S)/2$, only in cases where $S$ is compact,
		\item regions bounded by one or more geodesics or single neighboring curves in $S$, all of which are a single oriented 				distance from their respective isotopic geodesic in the direction of the least area set,
		\item horocusp neighborhoods when $S$ has cusps.
	\end{enumerate}

Adams and Morgan prove the existence and regularity of isoperimetric minimizers in all geometrically finite hyperbolic (constant section curvature $-1$) 2-manifolds \cite[Lemma 2.1]{AM}. Specifically, for each positive number $t$ less than the area of the manifold, they show the existence of a 2-dimensional isoperimetric minimizer of area $t$ whose boundary is comprised of embedded curves of equal constant curvature. Their proof uses the fact that if a 1-dimensional component goes off to infinity, either their enclosed area goes to 0, or it can be translated back into the interior of the surface. Using Corollary \ref{coro:Cheeger} we will prove the analogous result for Cheeger minimizers of a geometrically finite hyperbolic surface:
  \begin{theo}\label{theo:curvh}
   Let $S$ be a geometrically finite, finite area hyperbolic surface. Then there exists a 1-dimensional Cheeger minimizer in $S$ with positive length which is the boundary of a good Cheeger candidate of $S$.
  \end{theo}

A corollary of Theorem \ref{theo:curvh} is as follows:

\begin{coro}\label{coro:Surf}
		Let $S$ be a geometrically finite, finite area hyperbolic surface. Then the compact version of Buser's inequality and Theorem \ref{theo:AB} hold.
	\end{coro}

One application of Lemma \ref{lem:Stokes} is a proof of the known fact that the constant $\bar{h}$ of a horocusp neighborhood is equal to 1. Unfortunately, we do not know the original reference for this fact. We state this fact as a lemma and point out that this means that when $h(S)=1$, any horocusp neighborhood is a 2-dimensional Cheeger minimizer of $S$. We will then use this lemma to prove Theorem \ref{theo:curvh}.

\begin{lem}\label{lem:hcusp}
	Let $C$ be a horocusp neighborhood in $S$.  Then $\overline{h}(C)=h^{\ast}(C)=1$.
\end{lem}

\noindent {\bf Proof of Lemma \ref{lem:hcusp}.}  Let $\tilde{S}=\H^2$ denote the universal cover of $S$ where we denote the coordinates of $\tilde{S}$ to be the coordinates of the upper half-space $$\{(x,y) \in \R^2 : y>0\}$$ with the hyperbolic metric $ds^2=\frac{1}{y^2} \left (dx^2 + dy^2 \right )$.  We denote $\partial \tilde{S} = \{(x,y) \in \R^2: y=0\} \cup \{\infty\}$.

Via the isometries of $\H^2$, let $N$ be a neighborhood of cusp in $S$ where $\partial N$ pulls back to a horocircle around $\infty$ in $\tilde{S}\cup \partial \tilde{S}$ such that $N \subseteq C$.
Letting $p:\tilde{S} \ra S$ be the covering map, then we can parametrize coordinates of $p^{-1} (\partial N)$ by a curve $\gamma (t) = (t, 1)$ for $t\in [0,d]$ in the coordinates of the upper half-space.  Then we have 
	\begin{equation}\label{eq:areaN}
		\Area (N) = \int_{N} \frac{dA}{y^2} = \int_1^{\infty} \int_0^{d} \frac{dxdy}{y^2} = \int_1^{\infty}\frac{d}{y^2} dy= d
	\end{equation}
while
	\begin{equation}\label{eq:lengthN}
		\ell(\partial N) = \int_0^d \frac{\|\gamma'(t)\|}{\gamma_2 (t)} dt=\int_0^d ~dt = d.
	\end{equation}
Thus, $h^{\ast}(N) = \ell(\partial N)/\Area (N) = 1$ implying that $h(C) \leq 1$.

Let $D$ be any smooth subset of $C$. By showing that $h^{\ast} (D) \geq 1$ using Lemma \ref{lem:Stokes} with $a=1$, it will follow that $h^{\ast}(C) \geq 1$. Notice that if $\ell (\partial C) =L$, then we can parametrize $C$ on a ball of radius 1 in Euclidean polar coordinates $(r, \theta)$, with a metric of the form: $$ds^2=\frac{1}{r^2} dr^2+ \frac{L^2r^2}{4\pi^2} d\theta^2.$$ Then we can take $f_1(r,\theta)=1/r^2$ and $f_2(r,\theta)=L^2r^2/4\pi^2$. So then we have $f_1(r,\theta)f_2(r,\theta)=L/2\pi$. Taking $F(r,\theta)=f_2(r,\theta)$. Then by Lemma \ref{lem:Stokes}, we have that $h^{\ast}(D) \geq 1$. Since $D$ is an arbitrary smooth subset of $C$, we conclude that $\overline{h}(C)=1$.
\eproof

A very useful lemma for computing the region between a constant curvature curve and the neighboring geodesic in the same isotopy class is given by Adams and Morgan \cite{AM}:
	\begin{lem}\label{lem:AM} {\bf (Adams and Morgan \cite[Lemma 2.3]{AM})}
		Consider a curve of length $L$ and curvature $\kappa$ which neighbors a geodesic of length $L_0$ with a constant distance $s$ between them. Denote the area between the geodesic and the neighboring curve by $A$. Then we have
			$$L^2=A^2+L_0^2,$$
		and $A=L_0\sinh (s)$, $L=L_0 \cosh (s)$, and $\kappa = \tanh (s)$.
	\end{lem} 

Using Lemma \ref{lem:AM}, we can compute the isoperimetric ratio of an annulus. Then, by considering the isoperimetric ratio of the annulus with respect to the Cheeger constant, we prove the following:

	\begin{lem}\label{lem:restrict}
		Let $S$ be a geometrically finite, finite area hyperbolic surface. Then an annulus or the complement of an annulus in $S$ can only be a Cheeger minimizer when its area is exactly equal to $\Area(S)/2$ and $S$ is compact.		
	\end{lem}

\noindent {\bf Proof of Lemma \ref{lem:restrict}.}
We will argue that an annulus or the complement of an annulus can only be a Cheeger minimizer when their area is exactly equal to $\Area(S)/2$. By Lemma \ref{lem:AM}, the isoperimetric ratio of an annulus surrounding a geodesic of length $L_0$ is $\sqrt{1+ \left ( \frac{2L_0}{A}\right )^{\! 2}}$. Since this is an increasing function of $A$, we see that if a Cheeger minimizer of $S$ is an annulus, it would always seek to increase its area unless $A=\Area (S)/2$. For the complement of an annulus, we use Lemma \ref{lem:AM} to deduce that the isoperimetric ratio of the annulus complement is $$\sqrt{\left ( \frac{\Area(S)}{2A}-1 \right )^{\! 2} + \left ( \frac{2L_0}{A} \right )^{\! 2}}$$ where $A$ is the area of the annulus complement and $L_0$ is the length of the geodesic centered in the annulus. Since the isoperimetric ratio is a decreasing function of $A$, the stability of a Cheeger minimizer tells us that the annulus complement will always increase its area up to $A=\Area (S)/2$.

Finally, in the case where $S$ has a cusp, the isoperimetric ratio of a horocusp neighborhood is 1 which is strictly less than that of an annulus which is $\sqrt{1+ \left ( \frac{2L_0}{A}\right )^2}$, as noted above. Therefore an annulus cannot be a Cheeger minimizer when $S$ has a cusp. Since an annulus and its complement have the same isoperimetric when their areas are both equal to $\Area(S)/2$, this is also true for the complement of an annulus.\eproof

Denote the injectivity radius of a point $p \in S$ by $\inj_p (S)$. We write the injectivity radius of $S$ as $\inj (S)= \inf_{p \in S} \inj_p (S)$.

\begin{lem}\label{lem:restrict2}
	Let $S$ be a geometrically finite, finite area hyperbolic surface. If there is an embedded metric disk of radius $r$ centered at $p \in S$, denote it $D_r(p)$, so that $h^{\ast} \big (D_r(p) \big )=h(S)$ with $\Area \big (D_r(p) \big ) \leq \Area (S)/2$, then $\inj (S) \leq r$ and $\Area \big (D_r(p) \big ) =\Area (S)/2$.
\end{lem}

\noindent {\bf Proof of Lemma \ref{lem:restrict2}.} At any point $p \in S$, take $D_r(p)$ to be the metric ball centered at $p$ of radius $r$. Then the isoperimetric ratio of $D_r(p)$ is given by  $h^{\ast} \big (D_r(p)\big )=\frac{e^r-e^{-r}}{e^r+e^{-r}-2}$ and so $h^{\ast} \big (D_r(p)\big )$ is a decreasing function of radius $r$. Therefore, by the stability of a Cheeger minimizer, a metric ball $D_r(p)$ can only be a Cheeger minimizer when $\Area \big (D_r(p)\big ) = \Area(S)/2$. When $S$ is non-compact, we have that $h^{\ast}(C)=1$ for any horocusp neighborhood with $C \subset S$ with $\Area (C) \leq \Area (S)/2$. Further, since $h^{\ast} \big ( D_r(p)\big ) \searrow 1$ as $r \ra \infty$ and $\Area(S)<\infty$, we see that $D_r(p)$ cannot be a Cheeger minimizer of $S$.

Suppose that there is a point $p_1 \in S$ such that $D_r(p_1)$ is embedded but there exists $p_2 \in S$ such that $D_r(p_2)$ is not embedded. Take a smooth path $\gamma :[0,1] \ra S$ with $\gamma (0)=p_1$ and $\gamma (1)=p_2$. Then for $t \in [0,1]$, the injectivity radius $\inj_{\gamma (t)} (S)$ is a continuous function of $t$. By the intermediate value theorem, there exists $t' \in [0,1]$ such that the boundary $\partial D_r \big (\gamma(t') \big )$ is not embedded, but the interior of $D_r \big (\gamma(t') \big )$ is embedded. But this contradicts the embedded conclusion for boundary curves given in Theorem \ref{theo:Cheeger}.
\eproof

We are now ready to prove Theorem \ref{theo:curvh}.

\D {\bf Proof of Theorem \ref{theo:curvh}.} When $S$ is compact, then the existence and regularity portions of the result follow from Lemma \ref{lem:Buser} and Theorem \ref{theo:GMT3}.

When $S$ is non-compact, then there are two possibilities. If $h(S)=1$, then any horocusp neighborhood is a Cheeger minimizer, and the result follows. Now assume that $h(S)<1$. Then for each cusp $C_i \subset S$, consider a horocusp neighborhood $N_i$ of $C_i$. Then define $N= \bigsqcup_i N_i$. We have that $\overline{h}(N)= h^{\ast}(N)=1$ by Fact \ref{fact1} and Lemma \ref{lem:hcusp}. Since $\overline{h}(N)>h(S)$, $S$ is isoperimetrically compact with respect to the decomposition $(S-N) \cup N$. Then a Cheeger minimizer $A$ of $S$ exists by Corollary \ref{coro:Cheeger} and this Cheeger minimizer must have strictly positive area. It follows from the definition of the Cheeger constant that this Cheeger minimizer must also be a solution to the isoperimetric problem for $t=\Area (A)$.

\begin{center}
		\begin{figure}\label{fig:disk}
	\includegraphics[width=5in]{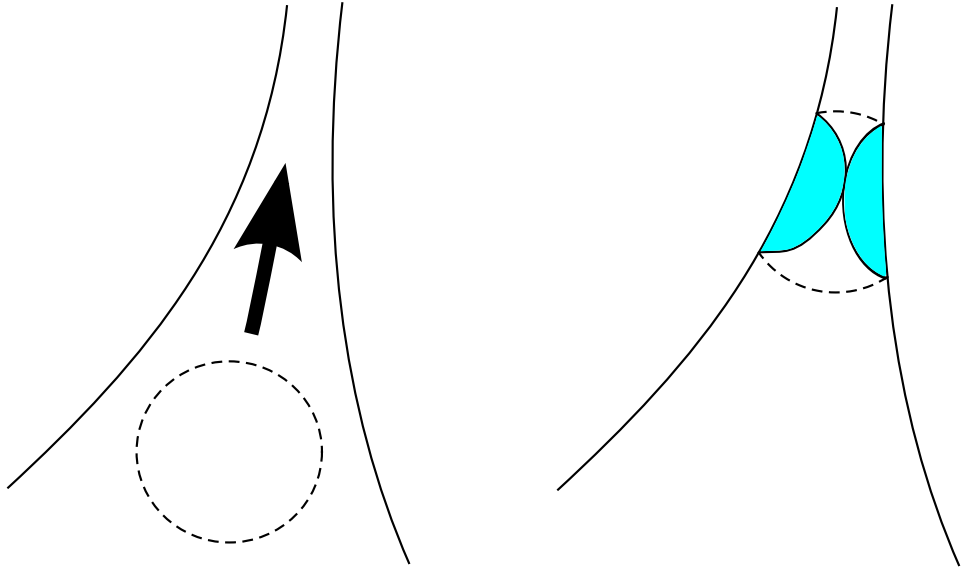}
        \caption{A disk positioned on the rear of the cusp is translated out of the cusp until its boundary meets at a single point.}\label{fig:CuspDisk}
		\end{figure}
	\end{center}

Lemma \ref{lem:restrict} tells us that an annulus (or its complement) can only be a Cheeger minimizer when its area is equal to $\Area (S)/2$. Further, when $S$ is compact, Lemma \ref{lem:restrict2} tells us that a disk can only be a Cheeger minimizer when its area is equal to $\Area (S)/2$. When $S$ is non-compact, any disk $D$ can be eliminated as a Cheeger candidate since its isoperimetric ratio is strictly greater than that of a horocusp neighborhood. The disk can also be eliminated because it cannot be embedded at every point in $S$, such as those in a horocusp neighborhood, see Figure \ref{fig:CuspDisk}, contradicting Theorem \ref{theo:Cheeger}. Therefore, a Cheeger minimizer must be a good Cheeger candidate.
\eproof

Adams and Morgan prove that there exists a hyperbolic metric on the genus two surface such that all isoperimetric minimizers are disks, annuli, complements of disks, or complements of annuli \cite[Theorem 3.1]{AM}. Since a Cheeger minimizer must also be an isoperimetric minimizer, it follows from Lemmas \ref{lem:restrict} and \ref{lem:restrict2} that their surface must either have a single  annulus or an embedded disk at every point having area equal to half the area of the surface. On the other hand, it is possible that there are very few examples of hyperbolic surfaces, if any, where a disk is a Cheeger minimizer. Consider the following reformulation of a conjecture about the systole of a closed hyperbolic surface: 
\begin{conj}\label{conj:SS}{\bf (Schmutz Schaller \cite[Conjecture 1]{SS98})}
	There exists a universal constant $C$ such that a surface $S_g$ of genus $g$ has an embedded metric disk $D_{\inj(S_g)}$ with
		\begin{equation}\label{eq:SS}
			\Area(D_{\inj(S_g)}) < 2\pi\left (3C\pi(g-1)^{2/3}-1 \right ).
		\end{equation}
\end{conj}

The constant $C$ in this conjecture is the same constant as is given in the original statement by Schmutz Schaller. We will show that if the conjecture above is true, there are only finitely many genera of hyperbolic surfaces which can contain an embedded disk whose area is half the area of the surface:
\begin{Claim}\label{Claim:disk}
	Assume that Conjecture \ref{conj:SS} is true. If $S_g$ is a hyperbolic surface of unspecified genus $g \geq 2$ with $C$ the universal constant in Conjecture \ref{conj:SS} and $\Area ( D_{\inj (S_g)} ) \geq \Area(S_g)/2$, then 
	$$g<27C^3+1.$$ 
\end{Claim}

\D {\bf Proof of Claim \ref{Claim:disk}.}
Suppose that $S$ is a hyperbolic surface with $\Area \left (D_{\inj(S)}\right ) \geq \Area (S)/2$. Then, by Gauss-Bonnet, we have that $\Area (S)/2 = 2\pi (g-1)$ and by equation (\ref{eq:SS}) from Conjecture \ref{conj:SS} we have that 
$$\Area \left (D_{\inj(S)}\right ) < 2\pi\left (3C(g-1)^{2/3}-1 \right )<2\pi \left (3C(g-1)^{\frac{2}{3}}\right ).$$ Combining these, we get that 
	\begin{equation}\label{eq:genus}
		g-1<3C(g-1)^{\frac{2}{3}}.
	\end{equation}
Solving equation (\ref{eq:genus}) for $g$, the result follows.
\eproof

%%%%%%%%%%%%%%%%%%%%%%%%%%%%%%%%%%%%%%%%%%%%%%%%%%%%%%%%%%%%%%%%%%%%%%%%%%%%%%%%%%%%
\section{An Algorithm for Computing the Cheeger Constant of $S$}
\label{sec:algorithm}
\setcounter{equation}{0}

Strohmaier and Uski gave an algorithm for computing the eigenvalues of the Laplacian on hyperbolic surfaces \cite{SU12}. Motivated by their result, we give an algorithm for computing the Cheeger constants of finite area hyperbolic surfaces. Adams and Morgan also use their characterization of isoperimetric minimizers of a geometrically finite hyperbolic surface $S$ to give an algorithm for finding the isoperimetric minimizers for the isoperimetric problem of a fixed area $t \in \big (0, \Area(S) \big )$ \cite[Section 4]{AM}. Because running their algorithm for all $t \in \big ( 0, \Area (S)/2 \big ]$ gives an uncountable number of total solutions corresponding to these isoperimetric problems, it does not directly translate to an algorithm to compute the Cheeger constant of $S$. We will show how to adapt their characterization of isoperiemtric minimizers in Theorem \ref{theo:GMT3} to give an algorithm for computing the Cheeger constant of a finite area hyperbolic surface $S$.

In order to show that the number of inputs for the algorithm are finite, we refer to the following well-known fact. 
\begin{lem}\label{lem:boundlen}
	Let $S$ be a geometrically finite, finite area hyperbolic surface. Then, for any $d>0$, there are finitely many simple closed geodesics $\gamma$ in $S$ such that $\ell (\gamma) \leq d$.
\end{lem}

We will sketch a proof of this lemma since we do not know of a reference to the literature for the case where $S$ is non-compact.

\noindent {\bf Proof of Lemma \ref{lem:boundlen}.} Buser provides a proof of the result for the case where $S$ is compact \cite[Theorem 1.6.11]{B92}. We will now assume that $S$ is non-compact.

Consider an ideal triangulation of $S$ where we denote $T_i$ to be the interior of the $i-th$ triangle. By the Gauss-Bonnet Theorem, no geodesic can be completely contained in a hyperbolic cusp. Let $\gamma$ be a geodesic in $S$ with $\ell (\gamma) \leq d$. Then there is a horocusp neighborhood $N$ contained in each cusp so that $\gamma \cap N =\emptyset$ for any such geodesic $\gamma$ of length less than $d$ by taking $N$ to have distance greater than $d$ from the base of the cusp. Since the vertices of $T_i$ are located at infinity in the cusps of $S$, its edges are geodesics, and its interior is simply-connected, there exists some $L_i>0$ corresponding to $T_i$ such that whenever $\gamma \cap \inter(T_i) \neq \emptyset$, we have that $\ell(\gamma \cap T_i) \geq L_i$.

The number of times that $\gamma$ crosses any triangle as a chord is bounded above by $d/\min_i\{L_i\}$. Since there are finitely many triangles, the number of times that $\gamma$ intersects an edge of any triangle is at most $2d/\min_i\{L_i\}$. So the number of edge intersection combinations that $\gamma$ can have is finite. Since the interior of each $T_i$ is simply-connected, the isotopy class of $\gamma$ is uniquely determined by the combinations of intersections that $\gamma$ makes with the edges of each triangle $T_i$. Since $\gamma$ is the unique geodesic in its isotopy class, the result follows.
\eproof

\subsection{Assumptions about computability}
For simplicity of presentation, we will write the algorithm as though we can compute with real numbers and we know the exact length of each geodesic curve and the distances between these curves in $S$. More formally, we could say that we are computing with a Blum-Shub-Smale machine; see Blum, Shub, and Smale for more details on computing with real numbers \cite{BSS}. In practice, however, we expect that any implementation of this algorithm will be constrained by the limitations of computing with a subset of the rational numbers and, possibly, by the approximations of the lengths of and distances between geodesic curves in $S$. The remainder of this subsection will suggest how this can be done from the presentation of our algorithm. 

Since one often cannot put the exact value of $\ell (\gamma_i)$ into a computer, one expects for the lengths of curves $\ell (\gamma_i)$ to be approximated by $L_i \in \Q_c$. Here $\Q_c$ is the set of rational numbers which can be represented in decimal form using the maximum number of floating point digits the computer and software can to perform the calculations prescribed by the algorithm. For each hyperbolic surface $S$, with input curves given by $\{\gamma_i\}$, we can think of approximating $\ell (\gamma_i)$ by $L_i \in \Q_c$. We want to specify that there is a very small $0<\epsilon \ll \min_i \ell(\gamma_i)$ so that $\ell(\gamma_i) \in (L_i-\epsilon, L_i +\epsilon)$. Specifically, we would like for $\Q_c$ to be a large enough subset of $\Q$ to allow for $\epsilon$ to be small enough to guarantee the following: whenever Cheeger candidates $A_1$ and $A_2$ are such that $h^{\ast} (A_1)<h^{\ast} (A_2)$, the isoperimetric ratios of their approximates also respect this ordering. If these assumptions hold, we can conclude that the region given by the algorithm to approximate $h(S)$ has a boundary consisting of curves in the same isotopy class as the curves in the boundary of an actual Cheeger minimizer.

\subsection{An algorithm for computing the Cheeger constant of a hyperbolic surface}
Consider a finite area hyperbolic surface $S$. We will assume that we know the area of $S$ and all combinations and lengths of embedded geodesics in $S$ whose total length summation is less than or equal to the supremum of $\ell (\partial A)$ for any Cheeger minimizer $A \subset S$. Note that by the length bound (\ref{eq:length}) in Theorem \ref{theo:GMT3} of Adams and Morgan and Theorem \ref{theo:curvh}, if $A$ is a Cheeger minimizer, then
	\begin{equation}\label{eq:lencom}
		\ell (\partial A) \leq \sqrt{\frac{\Area (S)^2}{4} + 2 \pi \Area (S)}.
	\end{equation} 
The following claim gives a sharper upper bound on $\ell (\partial A)$ when $S$ is non-compact:\begin{Claim}\label{Claim:length}
		 If $S$ is non-compact, then 
		 \begin{equation} \label{eq:lennon}
			\ell (\partial A) \leq \frac{\Area (S)}{2}.\end{equation}
\end{Claim}

\noindent {\bf Proof of Claim \ref{Claim:length}.}
 When $S$ is non-compact, then the isoperimetric ratio of a horocusp neighborhood provides an upper bound for $h^{\ast}(A)$. It follows that 
		$$\frac{\ell (\partial A)}{\Area (S)/2} 
		\leq \frac{\ell (\partial A)}{\Area (A)} \leq 1.$$
		This verifies Equation (\ref{eq:lennon}).		 
\eproof

Mirzakhani proved that given a fixed genus and number of cusps, that the number of geodesics of length less than or equal to $L$ grows polynomially in $L$ \cite[Theorem 1.1]{M08}. Therefore, we can conclude that the number of geodesics one needs to consider grows polynomially in $\Area(S)$.

We will enumerate these geodesics and denote them by $\gamma^i$ and think of them along with the following quantities as the inputs of our algorithm: we assume that we know the shortest perpendicular distance between each $\gamma^i$ and either another geodesic or itself in both normal directions from $\gamma^i$. We also assume that, for every splitting of $S$ into $A$ and $B$ via a collection of these geodesics, we know the Euler characteristics $\chi (A)$ and $\chi (B)$.

The algorithm outputs a real number $H$, which equals $h(S)$. It will also output a list of pairs of the form $\big ( \{\gamma^{i_1}, \ldots , \gamma^{i_j} \}, s_i \big )$ where $\gamma^{i_1}, \ldots , \gamma^{i_j}$ are geodesics whose union splits $S$ into $A$ and $B$ with $\Area (A) \leq \Area (B)$ and $s_i \geq 0$ is a distance from each $\gamma^{i_k}$ into $B$ which corresponds to a 
Cheeger minimizing hypersurface. In the case where $S$ has cusps and $H=h(S)=1$, one should note that any horocusp neighborhood is a 1-dimensional Cheeger minimizer.

\begin{enumerate}
	\item We will define initial values for $U$, the upper bound on the length of the sum of candidate geodesics, and $H$, the estimate for the Cheeger constant. If $S$ has a cusp, by Lemma \ref{lem:hcusp}, we know that the isoperimetric ratio of any horocusp neighborhood is 1, so set $H=1$ and $U=\Area (S)/2$ by Claim \ref{Claim:length}. If $S$ does not have any cusps, set $H= \infty$ and $U=\sqrt{\Area (S)^2/4+2\pi \Area (S)}$. Note that the latter choice of $U$ follows from Theorem \ref{theo:GMT3} of Adams and Morgan.
	\item Take a collection of geodesics $\gamma^{i_1}, \ldots , \gamma^{i_j}$ in $S$ which has not yet been considered by the algorithm. Specifically, this collection must have the following properties:
		\begin{enumerate}
			\item $\cup_{k=1}^j \gamma^{i_k}$ split $S$ into two subsurfaces $A$ and $B$ with 
			$A \cap B=\partial A = \partial B$ and
			\item $\sum_{k=1}^{j} \ell (\gamma^{i_k}) \leq U$.
		\end{enumerate}
		If all collections of geodesics satisfying these criteria have been considered, the algorithm terminates and returns the value of $H$ and pairs of the form $$\big ( \{\gamma^{i_1}, \ldots , \gamma^{i_j} \}, s_i \big ).$$
	\item When $\Vol_n(A)=\Vol_n(B)$, the isoperimetric ratio corresponding to the collection $\gamma^{i_1}, \ldots , \gamma^{i_j}$ is minimized and we can skip to step (5) where we can take $H_0=h^{\ast}(A)=h^{\ast} (B)$. Otherwise, proceed to step (4).
	\item Let $A_s$ and $B_s$ are the sets defined by the regions bounded by the constant curves of 			normal distance $s$ into $B$ where $A_0=A$ and $B_0=B$. It follows from the constant sectional curvature of $S$ and a first variation argument that the connected components of the boundary of an isoperimetric minimizer of $S$ must all be located at exactly the same distance to the unique closed geodesic in their respective isotopy classes. Determine the real number $d_{i_j}$ such that the level set of the signed distance function (positive into $B$) from the curve $\gamma^{i_j}$ at the distance of $d_{i_j}$ is not isotopic to the curve $\gamma^{i_j}$.
			By Lemma \ref{lem:AM} and the Gauss-Bonnet Theorem, we see that
		\begin{equation}\label{eq:exratios}
			h^{\ast}(A_s)=
				\frac{\ell (\partial A)\cosh(s)}{-2\pi \chi (A)+\ell (\partial A)\sinh (s)}, \quad
			h^{\ast} (B_s)= 
				\frac{\ell (\partial B)\cosh (s)}{-2\pi \chi (B) - \ell (\partial B)\sinh (s)},
		\end{equation}
	where $0 \leq 2s< \min_k d_{i_k}$. Since geodesics are length minimizing in their isotopy classes, we only need to consider sets where $A_s$ is increasing in area and $B_s$ is decreasing in area.
	\item Since $\ell (\partial B_s)$ is an increasing function of $s$, the ratio $h^{\ast}(B_s)$ can only obtain the minimum when $\Area (A_s)=\Area(B_s)$. Therefore, we need only minimize $h^{\ast}(A_s)$ over the single parameter $s \in \big [0, \min_k d_{i_k}/2 \big )$. It follows from (\ref{eq:exratios}) that the critical points of $h^{\ast}(A_s)$ as a function of $s$ are given when $\sinh(s) = -\frac{\ell (\partial A)}{2\pi \chi (A)}$. However, if the minimum is achieved at $\min_k d_{i_k}$, the region is not embedded and the corresponding candidate $A_s$ will not minimize $h^{\ast}(A_s)$ for all $s \in \big [0, \min_k d_{i_k}/2 \big )$. So in this case, we discard $A_s$ as a possible Cheeger minimizer. Now we consider the situation when minimum $H_0$ is achieved for some $s_0 \in \big [0, \min_k d_{i_k}/2 \big )$, which must happen for some collection $\gamma^{i_1}, \ldots , \gamma^{i_j}$ by Theorem \ref{theo:curvh}. If $H>H_0$ for this value of $s_0$, then redefine $H=H_0$ and replace all pairs with the single pair $\big ( \{i_1, \ldots , i_j \}, s_0 \big )$. If $H=H_0$, then add $\big ( \{i_1, \ldots , i_j \}, s_0\big )$ to the list of pairs. If $H < H_0$, do nothing.
	\item If $A$ is a 2-dimensional Cheeger minimizer, each value of $H$ that is obtained from step (5) of the algorithm gives us the inequality $H \geq \frac{2\ell (\partial A)}{\Area (S)}$. This gives the estimate $\ell (\partial A) \leq H \Area(S)/2$. Therefore, we compare this new estimate to our existing estimate for $U$. We conclude that if $H \Area (S)/2 <U$, we should redefine $U$ to be $H \Area (S)/2$; otherwise, we will leave $U$ unchanged. This is a standard branch and bound technique, which is well-known in the study of algorithms. In practice, this allows the algorithm to be more efficient than a brute-force algorithm.
	\item Return to step 2.
\end{enumerate}

After the algorithm terminates, it suffices to check the output against the isoperimetric ratios of a hyperbolic disk and annuli around each isotopy $\gamma^i$ all of which having area equal to $\Area(S)/2$. If any of these ratios are embedded and strictly less than the output $H$ of the algorithm, then redefine $H$ to be this ratio.\footnote{To simplify the presentation and because we do not know of any examples where this occurs, we do not include this as a distinct step in our algorithm.}

\begin{rmk} The Gauss-Bonnet Theorem for surfaces with cusps is the same as the version for compact surfaces. For a simple explanation of this fact, see Rosenberg \cite{SR82}.
\end{rmk}
%%%%%%%%%%%%%%%%%%%%%%%%%%%%%%%%%%%%%%%%%%%%%%%%%%%%%%%%%%%%%%%%%%%%%%%%%%%%%%%%%%%%%%%%%%%%%%%%%%%%%%%%%%%%%%%%%%%%%%%%%%%%%%%%%%%%%%%
\section{A Lower Bound on the Cheeger Constants of Arithmetic Surfaces}
\label{sec:Selberg}
\setcounter{equation}{0}

\subsection{Initial Assumptions} In this section, we will work under the assumption that $\lambda_1 \big (\omega (h) \big )$ in Theorem \ref{theo:AB} is a strictly increasing function of the Cheeger constant $h=h(S)$ of a surface $S$. Recall that $\omega (h)$ is a Sturm-Liouville problem which depends on a lower bound on the Ricci curvature of $S$ and $h(S)$. If we fix the lower bound on the Ricci curvature of $S$, we can think of $\omega (h)$ solely as a function of a real-valued parameter $h$, which corresponds to the Cheeger constant $h(S)$.

We will now define the Sturm-Liouville eigenvalue problem $\omega (h)$ for hyperbolic surfaces. First, consider an operator of the form
\begin{equation*}%\label{eq:SLop}
			L=-\frac{1}{J(\tau)}\frac{d}{d\tau} \left ( J(\tau) \frac{d}{d\tau} \right ),
\end{equation*}
where $J$ is a weight function defined by
	\begin{equation*}
		J(\tau)= \cosh(\tau) + h \sinh (\tau).
	\end{equation*}
Denote by $\omega(h)$ the SL problem on $(0,T)$ given by 
	\begin{equation*}
		Lu=\lambda u, \qquad u(0)=0, \qquad u'(T)=0.
	\end{equation*}
The endpoint $T$ is defined implicitly by the equation $$\frac{1}{h} = \int_0^T J(\tau) \, d\tau.$$

Our numerical experiments support the claim that $\lambda_1 \big (\omega (h) \big )$ is a strictly increasing function of $h \in \R_{>0}$, however, we currently do not have a rigorous proof. One difficulty is that as $h$ increases, the weight function $J$ increases while the right endpoint of the interval $T$, decreases. This makes a straight-forward comparison of the eigenvalues $\lambda_1 \big (\omega (h_1) \big )$ and $\lambda_1 \big ( \omega  (h_2) \big )$ for $h_2>h_1>0$ difficult.

\subsection{Selberg's Conjecture and Arithmetic Surfaces}
\label{sec:SelbergConj}

Let $\Gamma = \PSL (2, \Z)$, a discrete subgroup of $\PSL (2,\R)$.  Consider modular subgroups of the form
	\begin{equation}\label{eq:gamk}
		\Gamma_k := \{N \in \Gamma | N \equiv \pm \id \mod k\}
	\end{equation}
for $k \in \Z_{\geq 2}$, which are well-known to be torsion-free.  Brooks and Zuk consider hyperbolic surfaces of the form
	\begin{equation}\label{eq:Sk}
		S_k := \Gamma_k \backslash \H^2
	\end{equation}
for $k \in \Z_{\geq 2}$. Specifically, they prove the following result.
\begin{theo}\label{theo:BZ} {\bf (Brooks and Zuk \cite{BZ02})}  For all large $k \in \N$, we have that
	\begin{equation}\label{eq:BZhSk} h \left (S_k \right ) < 0.4402.
	\end{equation}
\end{theo}

An important conjecture was given by Selberg in 1965 and, as of this work, remains open \cite{S65}. For an introduction to Selberg's conjecture, see Sarnak \cite{S95}.
	\begin{conj}\label{conj:Selberg} {\bf (Selberg \cite{S65})}
		If $k \geq 1$, then $\lambda_1 (S_k) \geq \frac{1}{4}$.
	\end{conj}
Recent progress on this problem was made by Kim and Sarnak who give the following bound \cite{K03}:
	\begin{equation}\label{eq:KS}
		\lambda_1(S_k) \geq \frac{975}{4096}.
	\end{equation}
If we first assume that Selberg's conjecture true and apply Buser's inequality to $S_k$ with $\lambda_1=\frac{1}{4}$ giving a lower bound of $h(S_k) \geq 0.0707\ldots$, there is clearly no contradiction to Selberg's conjecture. Agol mentions in his concluding remarks that he would like to see this computation carried out for his quantitative improvement of Buser's inequality \cite{IA}.  Thus, we obtain the following result which contains the result of this calculation:

	\begin{theo}\label{theo:hS} Assume that $\lambda_1\big (\omega (h) \big )$ is monotone increasing as a function of $h \in \R^+$. Then $h(S_k)> 0.198727\ldots$. Assuming that Selberg's Conjecture is true, then $h(S_k) >0.205436\ldots$.
	\end{theo}

By comparing the second estimate in Theorem \ref{theo:hS} to Equation (\ref{eq:BZhSk}) of Theorem \ref{theo:BZ}, we find that these results are consistent with Selberg's conjecture.

\D {\bf Proof of Theorem \ref{theo:hS}.} First, equation (\ref{eq:KS}) of Kim and Sarnak combined with Theorem \ref{theo:AB} tell us that $\frac{975}{4096} \leq \lambda_1\big ( \omega (h) \big )$. We will use this equation to solve for a lower bound on $h(S_k)$. To do so, we consider the related initial value problem given by 
	\begin{equation} \label{eq:IVP1}
		u''(\tau)+\frac{J'(\tau)}{J(\tau)}u'(\tau)+\frac{975}{4096}u=0, \qquad u(T)=1, \qquad u'(T)=0.
	\end{equation}
Assuming that Selberg's conjecture is true, Theorem \ref{theo:AB} tells us that $\frac{1}{4} \leq \lambda_1\big ( \omega (h) \big )$. This gives the initial value problem
	\begin{equation} \label{eq:IVP2}
		u''(\tau)+\frac{J'(\tau)}{J(\tau)}u'(\tau)+\frac{1}{4}u=0, \qquad u(T)=1, \qquad u'(T)=0.
	\end{equation}
Note that the functions $u$ and $J$ as well as the endpoint $T$ all depend on $h(S)$. Therefore, we should think of solutions $u$ of (\ref{eq:IVP1}) as a one parameter family where the parameter corresponds to $h(S)$. Since (\ref{eq:IVP1}) is an initial value problem, existence and uniqueness of differential equations tells us that each choice of parameter $h$ gives a single unique solution $u$. Since $\omega(h)$ dictates that $u(0)=0$, the correct choice of $u$ is given by finding $h$ such that $u(0)=0$ is also true. The same observations are also true for the initial value problem given in (\ref{eq:IVP2}).\footnote{The solutions $u$ of (\ref{eq:IVP1}) and (\ref{eq:IVP2}) are not particularly useful for finding $h(M)$ numerically. See Benson for a detailed description of such solutions \cite{BB}.}

Using the SLEIGN2 program of Bailey, Everitt, and Zettl \cite{BEZ}, we determine that $h \approx  0.198727\ldots$ when $\lambda_1 (h) = \frac{975}{4096}$ and $h \approx 0.205436\ldots$ when $\lambda_1 (h)=\frac{1}{4}$. The conclusion follows from our assumption that $\lambda_1 \big (\omega (h) \big )$ is increasing as a real-valued function of $h \in \R^+$.
\eproof

%%%%%%%%%%%%%%%%%%%%%%%%%%%%%%%%%%%%%%%%%%%%%%%%%%%%%%%%%%%%%%%%%%%%%%%%%%%%%%%%%%%%%%%%%%%%%%%%%%%%%%%%%%%%%%%%%%%%%%%

\noindent {\bf Acknowledgements:} The author acknowledges support from National Science Foundation grant DMS 0838434 ``EMSW21MCTP: Research Experience for Graduate Students'' and grant DMS 1107452 ``RNMS: Geometric structures and representation varieties''. The author would also like to thank his advisor Nathan Dunfield for suggesting this project as well as for many helpful suggestions and insights. The author would also like to thank Ian Agol, Pierre Albin, Dave Auckly, Neal Coleman, David Dumas, Chris Judge, Grant Lakeland, Frank Morgan, Bruce Solomon, Brian White, and Luke Williams for helpful discussions. The author is very grateful to ICERM, where a portion of this work was completed while the author was a semester-long visitor.

\bibliographystyle{plain}
\bibliography{CheegerSurf}
\end{document}